# Quenched law of large numbers for branching Brownian motion in a random medium

János Engländer

a*Department of Statistics and Applied Probability, University of California, Santa Barbara CA 93106-3110, USA.
E-mail: englander@pstat.ucsb.edu*

**Abstract.** We study a spatial branching model, where the underlying motion is $d$-dimensional ($d \geq 1$) Brownian motion and the branching rate is affected by a random collection of reproduction suppressing sets dubbed *mild obstacles*. The main result of this paper is the quenched law of large numbers for the population for all $d \geq 1$. We also show that the branching Brownian motion with mild obstacles *spreads less quickly* than ordinary branching Brownian motion by giving an upper estimate on its speed. When the underlying motion is an arbitrary diffusion process, we obtain a *dichotomy* for the quenched local growth that is independent of the Poissonian intensity. More general offspring distributions (beyond the dyadic one considered in the main theorems) as well as mild obstacle models for superprocesses are also discussed.

**Résumé.** Nous étudions un modèle de branchement spatial où le mouvement de base est Brownien $d$-dimensionnel ($d \geq 1$) et le taux de branchement est modifié par une collection aléatoire d'ensembles sur lesquels la reproduction n'a pas lieu (*obstacles moux*). Le résultat principal de cet article est l'asymptotique (en probabilité) des taux de croissance globaux "quenchés" pour tout $d \geq 1$, et nous identifions les termes de correction sous-exponentielle. Nous montrons aussi que le branchement Brownien avec obstacles moux diffuse moins vite que le branchement Brownien classique en donnant une borne supérieure de sa vitesse. Dans le cas où le mouvement de base est un processus de diffusion arbitraire nous obtenons une *dichotomie* pour la croissance locale "quenchée" qui est indépendante de l'intensité Poissonienne. Le cas de distributions plus générales du nombre de descendants (autre que le cas dyadique considéré dans les théorème principaux), ainsi que des modèles d'obstacles moux pour des superprocessus, sont aussi discutés.

*MSC:* Primary 60J65; secondary 60J80; 60F10; 82B44

*Keywords:* Poissonian obstacles; Branching Brownian motion; Random environment; Fecundity selection; Radial speed; Wavefronts in random medium; Random KPP equation

# 1. Introduction

## 1.1. Model

The purpose of this article is to study a spatial branching model with the property that the branching rate is decreased in a certain random region. More specifically we will use a natural model for the random environment: *mild Poissonian obstacles*.

Let $\omega$ be a Poisson point process (PPP) on $\mathbb{R}^d$ with intensity $\nu > 0$ and let $\mathbf{P}$ denote the corresponding law. Furthermore, let $a > 0$ and $0 < \beta_1 < \beta_2$ be fixed. We define the *branching Brownian motion (BBM) with a mild Poissonian obstacle*, or the '$(\nu, \beta_1, \beta_2, a)$-BBM' as follows. Let $B(z, r)$ denote the open ball centered at $z \in \mathbb{R}^d$ with radius $r$ and let $K$ denote the random set given by the $a$-neighborhood of $\omega$:

$$K = K_\omega := \bigcup_{x_i \in \text{supp}(\omega)} \overline{B}(x_i, a).$$

Then $K$ is a *mild obstacle configuration* attached to $\omega$. This means that given $\omega$, we define $P^\omega$ as the law of the strictly dyadic (i.e. precisely two offspring) BBM on $\mathbb{R}^d$, $d \geq 1$, with spatially dependent branching rate

$$\beta(x, \omega) := \beta_1 1_{K_\omega}(x) + \beta_2 1_{K_\omega^c}(x).$$

An equivalent (informal) definition is that as long as a particle is in $K^c$, it obeys the branching rule with rate $\beta_2$, while in $K$ its reproduction is suppressed and it branches with the smaller rate $\beta_1$. (We assume that the process starts with a single particle at the origin.) The process under $P^\omega$ is called a *BBM with mild Poissonian obstacles* and denoted by $Z$. The total mass process will be denoted by $|Z|$. Further, $W$ will denote $d$-dimensional Brownian motion with probabilities $\{\mathbb{P}_x, \, x \in \mathbb{R}^d\}$.

As usual, for $0 < \gamma \leq 1$ and for non-empty compact sets $\mathcal{K} \subset \mathbb{R}^d$ define the set (Hölder-space) $C^\gamma(\mathcal{K})$, as the set of continuous bounded functions on $\mathcal{K}$ for which $\|f\|_{C^\gamma} := \|f\|_\infty + |f|_{C^\gamma}$ is finite, where $|f|_{C^\gamma} := \sup_{x \in \mathcal{K}, h \neq 0} \frac{|f(x+h) - f(x)|}{|h|^\gamma}$. Then $C^\gamma = C^\gamma(\mathbb{R}^d)$ will denote the space of functions on $\mathbb{R}^d$ which, restricted on $\mathcal{K}$, are in $C^\gamma(\mathcal{K})$ for all non-empty compact set $\mathcal{K} \subset \mathbb{R}^d$.

Finally, $\lfloor z \rfloor$ will denote the integer part of $z \in \mathbb{R}$: $\lfloor z \rfloor := \max\{n \in \mathbb{Z} \mid n \leq z\}$.

The following obvious *comparison* between an ordinary 'free' BBM and the one with mild obstacles may be useful to have in mind.

**Remark 1 (Comparison).** *Let $\widehat{P}$ correspond to the 'free' BBM with branching rate $\beta_2$ everywhere. Then for all $t \geq 0$, all $B \subseteq \mathbb{R}^d$ Borel, and all $k \in \mathbb{N}$,*

$$P^\omega(Z_t(B) < k) \geq \widehat{P}(Z_t(B) < k) \quad \mathbf{P}\text{-a.s.,}$$

*that is, the 'free' BBM is 'everywhere stochastically larger' than the BBM with mild obstacles, $\mathbf{P}$-a.s.*

## 1.2. Motivation

In this section we list some problems motivating the study of our particular setting.

(1) *Poissonian obstacles.* Topics concerning the *survival asymptotics* for a single Brownian particle *among Poissonian obstacles* became a classic subject in the last twenty years, initiated by Sznitman, Bolthausen and others. It was originally motivated by the celebrated 'Wiener-sausage asymptotics' of Donsker and Varadhan in 1975 and by mathematical physics, and the available literature today is huge – see the fundamental monograph [35] and the references therein; see also e.g. [3, 30, 34] for more recent results.

(2) *Spatial branching processes in a random environment.* In [11] a model of a *spatial branching process in a random environment* has been introduced. There, instead of mild obstacles, we introduced hard obstacles and instantaneous killing of the branching process once any particle hits the trap configuration $K$.

In [12] the model was further studied, and resulted in some unexpected results, and a similar model with 'individual killing rule' was also suggested. Individual killing rule means that only the individual particle hitting $K$ is eliminated and so the process dies when the last particle has been absorbed. The difference between that model and the model we consider in the present article is that now the killing mechanism is even 'milder' as it only suppresses temporarily the reproduction of the individual particle but does not eliminate the particle.

(3) *Law of large numbers.* It is a notoriously hard problem to prove the law of large numbers (that is, that the process behaves asymptotically as its expectation) for spatial branching systems – see e.g. [18] for a special case. In our situation the scaling is not purely exponential. In general, the not precisely exponential case is difficult[1] and it is interesting that in our setup *randomization* will help in proving a kind of LLN – see more in Section 3.1.

(4) *Wavefronts in random medium.* The spatial spread of our process is related to a work of Lee–Torcaso [31] on wave-front propagation for a random KPP equation and to earlier work of Freidlin [24] on the KPP equation with random coefficients – see Section 3.2 of this paper, and also [36] for a survey.

(5) *Catalytic spatial branching.* An alternative view on our setting is as follows. Arguably, the model can be viewed as a *catalytic* BBM as well – the catalytic set is then $K^c$ (in the sense that branching is 'intensified' there). Catalytic spatial branching (mostly for superprocesses though) has been the subject of vigorous research in the last twenty years initiated by Dawson, Fleischmann and others – see the survey papers [28] and [8] and references therein. In those models the individual branching rates of particles moving in space depend on the amount of contact between the particle ('reactant') and a certain random medium called the catalyst. The random medium is usually assumed to be a 'thin' random set (that could even be just one point) or another superprocess. Sometimes 'mutually' or even 'cyclically' catalytic branching is considered [8].

Our model is simpler than most catalytic models as our catalytic/blocking areas are fixed, whereas in several catalytic models they are moving. On the other hand, while for catalytic settings studied so far results were mostly only qualitative we are aiming to get a quite sharp *quantitative* result.

For the discrete setting there are much fewer results available. One example[2] is [26], where the branching particle system on $\mathbb{Z}^d$ is so that its branching is catalyzed by another autonomous particle system on $\mathbb{Z}^d$. There are two types of particles, the $A$-particles ('catalyst') and the $B$-particles ('reactant'). They move, branch and interact in the following way. Let $N_A(x,s)$ and $N_B(x,s)$ denote the number of $A$- [resp. $B$-]particles at $x \in \mathbb{Z}^d$ and at time $s \in [0,\infty)$. (All $N_A(x,0)$ and $N_B(x,0)$, $x \in \mathbb{Z}^d$ are independent Poisson variables with mean $\mu_A$ ($\mu_B$).)

Every $A$-particle ($B$-particle) performs independently a continuous-time random walk with jump rate $D_A$ ($D_B$). In addition a $B$-particle dies at rate $\delta$, and, when present at $x$ at time $s$, it splits into two in the next $\mathrm{d}s$ time with probability $\beta N_A(x,s)\,\mathrm{d}s + \mathrm{o}(\mathrm{d}s)$. Conditionally on the system of the $A$-particles, the jumps, deaths and splitting of the $B$-particles are independent. For large $\beta$ the existence of a critical $\delta$ is shown separating local extinction regime from local survival regime.

**Remark 2 (Self-duality).** *The discrete setting has the advantage that when $\mathbb{R}^d$ is replaced by $\mathbb{Z}^d$, the difference between the sets $K$ and $K^c$ is no longer relevant. Indeed, the equivalent of a Poisson trap configuration is an i.i.d. trap configuration on the lattice, and then its complement is also i.i.d. (with a different parameter). So, in the discrete case 'Poissonian mild obstacles' give the same type of model as 'Poissonian catalysts' would. This nice duality is lost in the continuous setting as the 'Swiss cheese' $K^c$ is not the same type of geometric object as $K$.*

(6) *Population models.* It appears that our proposed model of a BBM with 'mild' obstacles has biological merit to it (see Appendix B for more elaboration).

Returning to our mathematical model, consider the following natural questions (both in the annealed and the quenched sense):

(1) What can we say about the growth of the total population size?

---

[1]See for example the recent preprints [21] and [23].
[2]A further example of the discrete setting is [2].

(2) What are the large deviations? (e.g., what is the probability of producing an atypically small population?)
(3) What can we say about the *local* population growth?

As far as (1) is concerned, recall that the total population of an ordinary (free) BBM grows a.s. and in expectation as $e^{\beta_2 t}$. (Indeed, for ordinary BBM, the spatial component plays no role, and hence the total mass is just a $\beta_2$-rate pure birth process $X$. As is well known, the limit $N := \lim_{t \to \infty} e^{-\beta_2 t} X_t$ exists a.s. and in mean, and $P(0 < N < \infty) = 1$.) In our model of BBM with the reproduction blocking mechanism, how much will the suppressed branching in $K$ slow the global reproduction down? Will it actually change the exponent $\beta_2$? (We will see that although the global reproduction does slow down, the slowdown is captured by a sub-exponential factor, being different for the quenched and the annealed case.)

Consider now (2). Here is an argument to give a further motivation. Let us assume for simplicity that $\beta_1 = 0$ and ask the simplest question: what is the probability that there is no branching at all up to time $t > 0$? In order to avoid branching the first particle has to 'resist' the branching rate $\beta_2$ inside $K^c$. Therefore this question is quite similar to the survival asymptotics for a single Brownian motion among 'soft obstacles' – but of course in order to prevent branching the particle seeks for large islands covered by $K$ rather than the usual 'clearings.' In other words, the particle now prefers to avoid $K^c$ instead of $K$. Hence, (2) above is a possible generalization of this (modified) soft obstacle problem for a single particle, and the presence of branching creates new types of challenges.

As far as (3) is concerned, it will be shown to be related to the recent paper [15] which treats branching diffusions on Euclidean domains.

An additional problem, regarding the spatial *spread* of the process, will be investigated using the result obtained on the (quenched) growth of the total population size.

### *1.3. Outline*

The rest of this article is organized as follows. The next section gives an overview on some facts that are easy to derive from results in the literature. Section 3 presents the results, Theorems 1 and 2. Section 4 presents some further problems, while Section 5 contains some preparations and also gives the short proofs of the preliminary claims. In Section 6 the first main theorem is proven. Section 7 presents some additional problems – these problems are presented separately because they refer to certain parts of the proof. In Section 8 we prove Theorem 2. In the first appendix we explain how our mathematical setting relates to models in biology. Finally the second appendix presents a proof of a statement which is needed in the proofs – we think the statement is of independent interest.

## **2. Some preliminary claims**

In this section we state some results that are straightforward to derive from other results in the literature. The proofs are provided in Section 5.

### *2.1. Expected global growth and dichotomy for local growth*

Concerning the expected global growth rate we have the following result.

**Claim 3 (Expected global growth rate).** *On a set of full* **P**-*measure,*

$$E^\omega |Z_t| = \exp\left[\beta_2 t - c(d,\nu)\frac{t}{(\log t)^{2/d}}(1 + o(1))\right] \quad \text{as } t \to \infty \tag{1}$$

*(quenched asymptotics) and*

$$(\mathbf{E} \otimes E^\omega)|Z_t| = \exp[\beta_2 t - \tilde{c}(d,\nu)t^{d/(d+2)}(1 + o(1))] \quad \text{as } t \to \infty \tag{2}$$

*(annealed asymptotics), where*

$$c(d,\nu) := \lambda_d \left(\frac{d}{\nu\omega_d}\right)^{-2/d},$$

$$\tilde{c}(d,\nu) := (\nu\omega_d)^{2/(d+2)} \left(\frac{d+2}{2}\right) \left(\frac{2\lambda_d}{d}\right)^{d/(d+2)},$$

*and $\omega_d$ is the volume of the d-dimensional unit ball, while $\lambda_d$ is the principal Dirichlet eigenvalue of $-\frac{1}{2}\Delta$ on it.*

Notice that:

(1) $\beta_1$ does not appear in the formulas,
(2) the higher the dimension, the smaller the expected population size.

**Remark 4 (The difference between the quenched and the annealed asymptotics; dimension).** *Let us pretend for a moment that we are talking about an ordinary BBM. Then at time $t$ one has $e^{\beta_2 t}$ particles with probability tending to one as $t \to \infty$ (the population size divided by $e^{\beta_2 t}$ has a limit, to be precise). For $t$ fixed take a ball $B = B(0, R)$ (here $R = R(t)$) and let $K$ be so that $B \subset K^c$ (such a ball left empty by $K$ is called a* clearing*). Consider the expected number of particles that are confined to $B$ up to time $t$. These particles do not feel the blocking effect of $K$, while the other particles may have not been born due to it.*

*Optimize $R(t)$ with respect to the cost of having such a clearing and the probability of confining a single Brownian motion to it. This is precisely the same optimization as for the classical Wiener-sausage. Hence one gets the expectation in the theorem as a lower estimate.*

*One suspects that the main contribution in the expectation in (2) is coming from the expectation on the event of having a clearing with optimal radius $R(t)$. In other words, denoting by $p_t$ the probability that a single Brownian particle stays in the $R(t)$-ball up to time $t$, one argues heuristically that $p_t e^{\beta_2 t}$ particles will stay inside the clearing up to time $t$ 'for free' (i.e. with probability tending to one as $t \to \infty$).*

*The intuitive reasoning is as follows. If we had* independent *particles instead of BBM, then, by a 'Law of Large Numbers type argument' (using Chebyshev inequality and the fact that $\lim_{t \to \infty} p_t e^{\beta_2} = \infty$), roughly $p_t e^{\beta_2 t}$ particles out of the total $e^{\beta_2 t}$ would stay in the $R(t)$-ball up to time $t$ with probability tending to 1 as $t \uparrow \infty$. One suspects then that the lower estimate remains valid for the branching system too, because the particles are 'not too much correlated.' This kind of argument (in the quenched case though) will be made precise in the proof of our main theorem by estimating certain covariances.*

*To understand the difference between the annealed and the quenched case note that in the latter case large clearings (far away) are automatically (that is, $\mathbf{P}$-a.s.) present. Hence, similarly to the single Brownian particle problem, the difference between the two asymptotics is due to the fact that even though there is an appropriate clearing far away $\mathbf{P}$-a.s., there is one around the origin with a small (but not too small) probability. Still, the two cases will have a similar element when, in Theorem 1 we drop the expectations and investigate the process itself, and show that inside such a clearing a large population is going to flourish (see the proof of Theorem 1 for more on this). This also leads to the intuitive explanation for the* decrease *of the population size as the dimension increases. The radial drift in dimension $d$ is $(d-1)/2$. The more transient the motion (i.e., the larger $d$), the harder it is for particles to stay in the appropriate clearings.*

Concerning local population size we present the following (quenched) claim.

**Claim 5 (Quenched exponential growth).** *The following holds on a set of full $\mathbf{P}$-measure: For any $\varepsilon > 0$ and any bounded open set $\emptyset \neq B \subset \mathbb{R}^d$,*

$$P^\omega_\mu \left( \limsup_{t \uparrow \infty} e^{-(\beta_2 - \varepsilon)t} Z_t(B) = \infty \right) > 0 \quad \text{and} \quad P^\omega_\mu \left( \limsup_{t \uparrow \infty} e^{-\beta_2 t} Z_t(B) < \infty \right) = 1.$$

**Problem 6.** What can one say about the distribution of the global and the local population size? The quenched asymptotics of the expected logarithmic global growth rate suggests that perhaps the rate itself is around $\beta_2 - c(d,\nu)(\log t)^{-2/d}$ at time $t$ in some sense. This problem will be addressed in Section 3.1. We will prove an appropriate formulation of the statement when the limit is meant in probability.

We now show how Claim 5 can be generalized for the case when the underlying motion is a diffusion. Let **P** be as before but replace the Brownian motion by an $L$-diffusion $X$ on $\mathbb{R}^d$, where $L$ is a second-order elliptic operator of the form

$$L = \sum_{i,j=1}^n a_{i,j}(x)\frac{\partial^2}{\partial x_i \partial x_j} + \sum_{i=1}^n b_i(x)\frac{\partial}{\partial x_i}$$

with $a_{i,j}, b_i \in C^\gamma$, $\gamma \in [0,1)$, and the symmetric matrix $a_{i,j}(x)$ is positive definite for all $x \in \mathbb{R}^d$. The branching $L$-diffusion with the Poissonian obstacles can be defined analogously to the case of BM. To present the result, we need an additional concept. Let

$$\lambda_c(L) = \lambda_c(L, \mathbb{R}^d) := \inf\{\lambda \in \mathbb{R} : \exists u > 0 \text{ satisfying } (L-\lambda)u = 0 \text{ in } \mathbb{R}^d\}$$

denote the *generalized principal eigenvalue* for $L$ on $\mathbb{R}^d$. In fact $\lambda_c \leq 0$, because $L1 = 0$ – see Section 4.4 in [32].

The following claim shows that the local behavior of the process exhibits a *dichotomy*. The crossover is given in terms of the local branching rate $\beta_2$ and the term $\lambda_c(L)$: local extinction occurs when the branching rate inside the 'free region' $K^c$ is not sufficiently large to compensate the transience of the underlying $L$-diffusion; if it is strong enough, then local mass grows exponentially. Note an interesting feature of the result: *neither $\beta_1$ nor the intensity $\nu$ of the obstacles play a role.*

**Claim 7 (Quenched exponential growth/local extinction).** *Given the environment $\omega$, denote by $P^\omega$ the law of the branching $L$-diffusion.*

(i) *Let $\beta_2 > -\lambda_c(L)$ and let $\nu > 0$ and $\beta_1 \in (0, \beta_2)$ be arbitrary. Then the following holds on a set of full **P**-measure: For any $\varepsilon > 0$ and any bounded open set $\emptyset \neq B \subset \mathbb{R}^d$,*

$$P^\omega\left(\limsup_{t\uparrow\infty} e^{(-\beta_2-\lambda_c(L)+\varepsilon)t} Z_t(B) = \infty\right) > 0$$

*and*

$$P^\omega\left(\limsup_{t\uparrow\infty} e^{(-\beta_2-\lambda_c(L))t} Z_t(B) < \infty\right) = 1.$$

(ii) *Let $\beta_2 \leq -\lambda_c(L)$ and let $\nu > 0$ and $\beta_1 \in (0, \beta_2)$ be arbitrary. Then the following holds on a set of full **P**-measure: For any bounded open set $B \subset \mathbb{R}^d$ there exists a $P^\omega$-a.s. finite random time $t_0 = t_0(B)$ such that $X_t(B) = 0$ for all $t \geq t_0$ (local extinction).*

## 3. Results

*3.1. Quenched asymptotics of global growth; LLN*

In this section we investigate the behavior of the (quenched) global growth rate.

As already mentioned in the Introduction of this paper, it is a notoriously hard problem to prove the law of large numbers for spatial branching systems, and, in general, the not purely exponential case is very challenging.

To elucidate this point, let us consider a branching diffusion with motion component corresponding to a second order elliptic operator $L$ on some Euclidean domain $D \subset \mathbb{R}^d$ and with spatially dependent (not identically zero) branching rate $\beta \geq 0$. Let $0 \not\equiv f$ be a nonnegative compactly supported bounded measurable function on $D$. If $\lambda_c$, the generalized principal eigenvalue of $L + \beta$ on $D$ is positive and $S = \{S_t\}_{t \geq 0}$ denotes the semigroup corresponding to $L + \beta$ on $D$, then $(S_t f)(\cdot)$ grows (pointwise) as $e^{\lambda_c t}$ on an exponential scale. However, in general, the scaling is not precisely exponential due to the presence of a subexponential term.

Instead of discussing here the appropriate spectral conditions for having precisely exponential scaling we refer the reader to Appendix A of [19] and note that in the present case it turns out that the resulting operator $\frac{1}{2}\Delta + \beta$ does not scale precisely exponentially **P**-a.s. Replacing $f$ by the function $g \equiv 1$, it is still true that the growth is not precisely exponential – this is readily seen in Claim 3 and its proof.

Since the process in expectation is related to the semigroup $S$, therefore purely exponential scaling indicates the same type of scaling for the expectation of the process (locally or globally). It turns out that if one is interested in the scaling of the process itself (not just the expectation), the case when there is an additional subexponential factor is much harder (see e.g. [19] and references therein).

In the proof of Theorem 1 below it is the *randomization* of the branching rate $\beta$ that helps, as $\beta$ has some 'nice' properties for almost every environment $\omega$, i.e. the 'irregular' branching rates 'sit in the **P**-zero set.'

Define the *average growth rate* by

$$r_t = r_t(\omega) := \frac{\log |Z_t(\omega)|}{t}.$$

Replace now $|Z_t(\omega)|$ by its expectation $\overline{Z}_t := E^\omega |Z_t(\omega)|$ and define

$$\widehat{r}_t = \widehat{r}_t(\omega) := \frac{\log \overline{Z}_t}{t}.$$

Recall from Claim 1, that on a set of full **P**-measure,

$$\lim_{t \to \infty} (\log t)^{2/d} (\widehat{r}_t - \beta_2) = -c(d, \nu). \tag{3}$$

We are going to show that an analogous statement holds for $r_t$ itself.

**Theorem 1 (LLN).** *On a set of full **P**-measure,*

$$\lim_{t \to \infty} (\log t)^{2/d} (r_t - \beta_2) = -c(d, \nu) \quad \text{in } P^\omega\text{-probability.} \tag{4}$$

One interprets Theorem 1 as a kind of quenched law of large numbers. Loosely speaking,

$$r_t \approx \beta_2 - c(d, \nu)(\log t)^{-2/d} \approx \widehat{r}_t, \quad t \to \infty,$$

on a set of full **P**-measure.

*3.2. The spatial spread of the process*

*3.2.1. Introduction*

L. Mytnik asked *how much the speed (spatial spread) of free BBM reduces* due to the presence of the mild obstacle configuration. Note that we are not talking about the bulk of the population (or the 'shape') but rather about *individual* particles traveling to very large distances from the origin (cf. Section 7).

As is well known, ordinary 'free' branching Brownian motion with constant branching rate $\beta_2 > 0$ has radial speed $\sqrt{2\beta_2}$ (see the recent paper [29] for a review and a strong version of this statement). Let $N_t$ denote the population size at $t \geq 0$ and let $\xi_k$ ($1 \leq k \leq N_t$) denote the position of the $k$th particle (with arbitrary labeling) in the population. Furthermore, let $m(t)$ denote a number for which $u(t, m(t)) = \frac{1}{2}$, where

$$u(t, x) := P\Big[\max_{1 \leq k \leq N_t} \|\xi_k(t)\| \leq x\Big].$$

In his classic paper, Bramson [5] considered the one-dimensional case and proved that as $t \to \infty$,

$$m(t) = t\sqrt{2\beta_2} - \frac{3}{2\sqrt{2\beta_2}} \log t + \mathcal{O}(1). \tag{5}$$

Since the one-dimensional projection of a $d$-dimensional branching Brownian motion is a one-dimensional branching Brownian motion, we can utilize Bramson's result for the higher dimensional cases too. Namely, it is clear, that in high dimension the spread is *at least* as quick as in (5). In [5] the asymptotics (5) is derived for the case $\beta_2 = 1$; the general result can be obtained similarly. See also p. 438 in [24]. (It is also interesting to take a look at [4].) Studying the function $u$ has significance in analysis too as $u$ solves

$$\frac{\partial u}{\partial t} = \frac{1}{2} u_{xx} + \beta_2(u^2 - u), \tag{6}$$

with initial condition

$$u(0, \cdot) = 1_{[0,\infty)}(\cdot). \tag{7}$$

In this section we show that the branching Brownian motion with mild obstacles *spreads less quickly* than ordinary branching Brownian motion by giving an upper estimate on its speed.

A related result was obtained earlier by Lee and Torcaso [31], but, unlike in (5), only up to the linear term and moreover, for random walks instead of Brownian motions. The approach in [31] is to consider the problem as the description of wavefront propagation for a random KPP equation. They extended a result of Freidlin and Gärtner for KPP wave fronts to the case $d \geq 2$ for i.i.d. random media. In [31] the wavefront propagation speed is attained for the discrete-space (lattice) KPP using a large deviation approach. Note that the 'speed' is only defined in a logarithmic sense. More precisely, let $u$ denote the solution of the discrete-space KPP equation with an initial condition that vanishes everywhere except the origin. The authors define a bivariate function $F$ on $\mathbb{R} \times \mathbb{R}^d \setminus \{\mathbf{0}\}$ and show that it satisfies

$$\lim_{t \to \infty} \frac{1}{t} \log u(t, tv\mathbf{e}) = -[F(v, \mathbf{e}) \vee 0],$$

for all $v > 0$ and $\mathbf{e} \in \mathbb{R}^d \setminus \{\mathbf{0}\}$. It turns out that there is a unique solution $v = v_{\mathbf{e}}$ to $F(v, \mathbf{e}) = 0$, and $v_{\mathbf{e}}$ defines the 'wave speed.' In particular, the speed is non-random.

Unfortunately, it does not seem to be easy to evaluate the variational formula for $v_{\mathbf{e}}$ given in [31], even in very simple cases.

It should be pointed out that the problem is greatly simplified for $d = 1$ and it had already been investigated by Freidlin in his classic text [24] (the *KPP equation with random coefficients* is treated in Section VII.7.7.). Again, it does not seem clear whether one can easily extract an explicit result for the speed of a branching RW with i.i.d. branching coefficients which can only take two values, $0 < \beta_1 < \beta_2$ (bistable nonlinear term).

The description of wavefronts in random medium for $d > 1$ is still an open and very interesting problem. The above work of Torcaso and Lee concerning processes on $\mathbb{Z}^d$ is the only relevant article we are aware of. To the best of our knowledge, the problem is open; it is of special interest for a bistable nonlinear term.

Before turning to the upper estimate, we discuss the lower estimate.

*3.2.2. On the lower estimate on the radial speed*
We are going to show that, if in our model Brownian motion is replaced by Brownian motion with constant drift $\gamma$ in a given direction, then any fixed non-empty ball is recharged infinitely often with positive probability, as long as the drift satisfies $|\gamma| < \sqrt{2\beta_2}$.

For simplicity, assume that $d = 1$ (the general case is similar). Fix the environment $\omega$. Recall Doob's $h$-transform of second-order elliptic operators:

$$L^h(\cdot) := \frac{1}{h} L(h \cdot).$$

Applying an $h$-transform with $h(x) := \exp(-\gamma x)$, a straightforward computation shows that the operator

$$L := \frac{1}{2}\frac{\mathrm{d}^2}{\mathrm{d}x^2} + \gamma\frac{\mathrm{d}}{\mathrm{d}x} + \beta_2$$

transforms into

$$L^h = \frac{1}{2}\frac{\mathrm{d}^2}{\mathrm{d}x^2} - \frac{\gamma^2}{2} + \beta_2.$$

Then, similarly to the proof of Claim 7, one can show that the generalized principal eigenvalue for this latter operator is $-\frac{\gamma^2}{2} + \beta_2$ for almost every environment. Since the generalized principal eigenvalue is invariant under $h$-transforms, it follows that $-\frac{\gamma^2}{2} + \beta_2 > 0$ is the generalized principal eigenvalue of $L$. Hence, by Lemma 10, any fixed non-empty interval is recharged infinitely often with positive probability.

Turning back to our original setting, the application of the '*spine*'-technology seems also promising.

Roughly speaking, the spine method uses a Girsanov-type change of measure on the branching diffusion. On one hand, the Girsanov density is a known (mean-one) martingale, for which the question of uniform integrability can be checked. On the other hand, under the new measure, the process can be represented 'directly' by a 'spine' or 'backbone' type construction. The 'spine' represents a special '*deviant*' particle performing a motion which has important properties from the point of view of the problem regarding the original branching process. The question regarding 'going back to the original measure' is then related to the UI property of the martingale. (See [25, 29] and references therein.)

In our case it is probably not very difficult to show the existence of a 'spine' particle (under a martingale-change of measure) that has drift $\gamma$ as long as $|\gamma| < \sqrt{2\beta_2}$.

### 3.2.3. An upper estimate on the radial speed

Our main result in this section is an upper estimate on the speed of the process. We give an upper estimate in which the order of the correction term is *larger* than the $\mathcal{O}(\log t)$ term appearing in Bramson's result, namely it is $\mathcal{O}(\frac{t}{(\log t)^{2/d}})$. (All orders are meant for $t \to \infty$.) We show that, loosely speaking, at time $t$ the spread of the process is not more than

$$t\sqrt{2\beta_2} - c(d,\nu)\sqrt{\frac{\beta_2}{2}} \cdot \frac{t}{(\log t)^{2/d}}.$$

(Again, $\beta_1$ plays no role as long as $\beta_1 \in (0, \beta_2)$.) The precise statement is as follows.

**Theorem 2.** *Define the functions*

$$f(t) := c(d,\nu)\frac{t}{(\log t)^{2/d}} \quad and \quad n(t) := t\sqrt{2\beta_2} \cdot \sqrt{1 - \frac{f(t)}{\beta_2 t}},$$

*where we recall that $c(d,\nu) := \lambda_d(\frac{\nu\omega_d}{d})^{2/d}$ and $\omega_d$ is the volume of the $d$-dimensional unit ball, while $\lambda_d$ is the principal Dirichlet eigenvalue of $-\frac{1}{2}\Delta$ on it. Then*

$$n(t) = t\sqrt{2\beta_2} - c(d,\nu)\sqrt{\frac{\beta_2}{2}} \cdot \frac{t}{(\log t)^{2/d}} + \mathcal{O}\left(\frac{t}{(\log t)^{4/d}}\right). \tag{8}$$

*Furthermore, if*

$$A_t := \{\text{no particle has left the } n(t) - \text{ball up to } t\}$$
$$= \left\{\bigcup_{0 \leq s \leq t} \mathrm{supp}(Z_s) \subseteq B(0, n(t))\right\},$$

*then on a set of full* **P**-*measure,*

$$\liminf_{t \to \infty} P^\omega(A_t) > 0. \tag{9}$$

## 4. Further problems

In this section we suggest some further problems and directions for research.

### 4.1. More general branching

It should also be investigated, what happens when the dyadic branching law is replaced by a general one (but the random branching rate is as before). In a more sophisticated population model, particles can also die – then the obstacles do not necessarily reduce the population size as they sometimes prevent death.

#### 4.1.1. Supercritical branching
When the offspring distribution is supercritical, the method of our paper seems to work, although when the offspring number can also be zero, one has to *condition on survival* for getting the asymptotic behavior.

#### 4.1.2. (Sub)critical branching
Critical branching requires an approach very different from the supercritical one, since taking expectations now does not provide a clue: $E^\omega |Z_t(\omega)| = 1$, $\forall t > 0$, $\forall \omega \in \Omega$.

Having the obstacles, the first question is whether it is still true that

$$P^\omega(\text{extinction}) = 1 \quad \forall \omega \in \Omega.$$

The answer is yes. To see this, note that since $|Z|$ is still a martingale, it has a nonnegative a.s. limit. This limit must be zero; otherwise $|Z|$ would stabilize at a positive integer. This, however, is impossible because following one Brownian particle it is obvious that this particle experiences branching events for arbitrarily large times.

Setting $\beta_1 = 0$, the previous argument still goes through. Let $\tau$ denote the almost surely finite extinction time for this case. One of the basic questions is the decay rate for $P^\omega(\tau > t)$. Will the tail be significantly heavier than $\mathcal{O}(1/t)$? (Of course $1/t$ would be the rate without obstacles.)

The subcritical case can be treated in a similar fashion. In particular, the total mass is a supermartingale and $P^\omega(\text{extinction}) = 1$ $\forall \omega \in \Omega$.

### 4.2. Superprocesses with mild obstacles

A further goal is to generalize the setting by defining *superprocesses with mild obstacles* analogously to the BBM with mild obstacles. Recall the concept of an $(L, \beta, \alpha)$-superdiffusion: Let $\mathcal{M}_f$ denote the set of finite measures $\mu$ on $\mathbb{R}^d$ and let $\alpha, \beta$ denote functions in the Hölder space $C^\gamma$ satisfying $\alpha > 0$ and $\sup_{x \in \mathbb{R}^d} \beta(x) < \infty$.

**Notation 8.** *Let $(X, \mathbf{P}_\mu, \mu \in \mathcal{M}_f)$ denote the $(L, \beta, \alpha)$-superdiffusion. That is, $X$ is the unique $\mathcal{M}_f$-valued (time-homogeneous) continuous Markov process which satisfies, for any bounded continuous $g : \mathbb{R}^d \mapsto \mathbb{R}_+$,*

$$\mathbf{E}_\mu \exp\langle X_t, -g\rangle = \exp\left(-\int_{\mathbb{R}^d} \mu(\mathrm{d}x)\, u(t, x)\right), \tag{10}$$

*where $u$ is the minimal nonnegative solution to*

$$\left. \begin{array}{l} \frac{\partial}{\partial t} u = Lu + \beta u - \alpha u^2 \quad \text{on } \mathbb{R}^d \times (0, \infty), \\ \lim_{t \to 0+} u(t, \cdot) = g(\cdot). \end{array} \right\} \tag{11}$$

*(Here $\langle \nu, f\rangle$ denotes the integral $\int_{\mathbb{R}^d} \nu(\mathrm{d}x)\, f(x)$.)*

One usually refers to $L$ as *migration*, $\beta$ as *mass creation* and $\alpha$ as the *intensity parameter* (or variance).

In fact, $X$ also arises as the short life time and high density diffusion limit of a *branching particle system*, which can be described as follows: in the $n$th approximation step each particle has mass $1/n$ and lives a random time which is exponential with mean $1/n$. While a particle is alive, its motion is described by a diffusion process corresponding to the operator $L$. At the end of its life, the particle dies and is replaced by a random number of particles situated at the parent particle's final position. The distribution of the number of descendants is spatially varying such that the mean number of descendants is $1 + \frac{\beta(x)}{n}$, while the variance is assumed to be $2\alpha(x)$. All these mechanisms are independent of each other.

The definition of the superprocess with mild obstacles is straightforward: the parameter $\alpha$ on the (random) set $K$ is smaller than elsewhere.

Similarly, one can consider the case when instead of $\alpha$, the 'mass creation term' $\beta$ is random, for example with $\beta$ defined in the same way (or with a mollified version) as for the discrete branching particle system. Denote now by $P^\omega$ the law of this latter superprocess for a given environment. We suspect that the superprocess with mild obstacles behaves similarly to the discrete branching process with mild obstacles when $\lambda_c(L+\beta) > 0$ and $P^\omega(\cdot)$ is replaced by $P^\omega(\cdot | X \text{ survives})$. The upper estimate can be carried out in a manner similar to the discrete particle system, as the expectation formula is still in force for superprocesses.

As we have already pointed out, there is a large amount of ongoing research on catalytic superprocesses; $\alpha$ is usually taken as a thin (sometimes randomly moving) set, or even another superprocess. In those models, one usually cannot derive sharp quantitative results. In a very simple one-dimensional model, introduced in [13], $\beta$ was spatially varying but deterministic and non-moving – in fact it was the *Dirac delta at zero*. Nevertheless, already in this simple model it was quite challenging to prove the asymptotic behavior of the process (Theorem 2 in [18]). Fleischmann, Mueller and Vogt suggest, as an open problem, the description of the asymptotic behavior of the process in the *three-dimensional* case [23]; the two-dimensional case is even harder, as getting the asymptotics of the *expectation* is already difficult. As we have mentioned, the *randomization* of $\beta$ may help in the sense that $\beta$ has some 'nice' properties for almost every environment $\omega$.

## 5. Preparations

In this section we present some preparatory lemmas and also the proofs of the preliminary claims.

**Lemma 9 (Expectation given by Brownian functional).** *Fix $\omega$. Then*

$$E^\omega |Z_t| = \mathbb{E} \exp\left[\int_0^t \beta(W_s)\right] ds. \qquad (12)$$

**Proof.** It is well known ('first moment formula' of spatial branching processes[3]; see, for instance (1.50a) in [10]), that $E_x^\omega |Z_t| = (T_t 1)(x)$, where $u(x,t) := (T_t 1)(x)$ is the minimal solution of the parabolic problem:

$$\frac{\partial u}{\partial t} = \left(\frac{1}{2}\Delta + \beta\right) u \quad \text{on } \mathbb{R}^d \times (0, \infty),$$
$$u(\cdot, 0) = 1, \qquad (13)$$
$$u \geq 0.$$

(Here $\{T_t\}_{t \geq 0}$ denotes of course the semigroup corresponding to the generator $\frac{1}{2}\Delta + \beta$ on $\mathbb{R}^d$.) This is equivalent (by the Feynman–Kac formula) to (12). □

We will also need the following result: Let $D \subseteq \mathbb{R}^d$ be a (non-empty) domain and let $0 \leq V$, $0 \not\equiv V$ be bounded from above and in $C^\gamma(D)$, $\gamma \in (0,1]$. The (binary) $(L, V)$-*branching diffusion* (or more precisely, the $(L, V; D)$-branching diffusion) is the Markov process with motion component $Y$ and with spatially dependent branching rate $V$, replacing particles by precisely two offspring when branching. At each time $t > 0$, the process consists of a point process $X_t$ defined on Borel sets of $D$.

---

[3]The statements and their derivations are the same for branching processes and superprocesses.

**Lemma 10 ([15], Theorem 3).** *Given $D$, let $P$ denote the law of the $(L, V)$-branching diffusion $X$ and let $\lambda_c := \lambda_c(L + V)$.*

(i) *Under $P$ the process $X$ exhibits local extinction if and only if $\lambda_c \leq 0$.*
(ii) *When $\lambda_c > 0$, for any $\lambda < \lambda_c$ and any open $\emptyset \neq B \subset\subset D$, one has*[4]

$$P_\mu\Big(\limsup_{t\uparrow\infty} e^{-\lambda t} X_t(B) = \infty\Big) > 0 \quad \text{and} \quad P_\mu\Big(\limsup_{t\uparrow\infty} e^{-\lambda_c t} X_t(B) < \infty\Big) = 1.$$

*5.1. Proof of Lemma 9*

Since $\beta := \beta_1 1_K + \beta_2 1_{K^c} = \beta_2 - (\beta_2 - \beta_1) 1_K$, we can rewrite the equation (12) as

$$E^\omega |Z_t| = e^{\beta_2 t} \mathbb{E} \exp\Big[-\int_0^t (\beta_2 - \beta_1) 1_K(W_s) \, ds\Big].$$

The expectation on the right-hand side is precisely the survival probability among soft obstacles of 'height' $\beta_2 - \beta_1$, except that we do not sum the shape functions on the overlapping balls. However this does not make any difference with regard to the asymptotics (see [35], Remark 4.2.2.). The statements thus follow from the well known Donsker–Varadhan type 'Wiener asymptotics' for soft obstacles (see [35], Theorems 4.5.1. and 4.5.3.; see also [9]).

*5.2. Proof of Claim 5*

Since in this case $L = \Delta/2$ and since $\lambda_c(\Delta/2, \mathbb{R}^d) = 0$, the statement immediately follows from Claim 7.

*5.3. Proof of Claim 7*

In order to be able to use Lemma 10, we compare the rate $\beta$ with another, *smooth* (i.e. $C^\gamma$) function $V$. Recalling that $K = K_\omega := \bigcup_{x_i \in \text{supp}(\omega)} \overline{B}(x_i, a)$, let us enlarge the obstacles:

$$K^* = K_\omega^* := \bigcup_{x_i \in \text{supp}(\omega)} \overline{B}(x_i, 2a).$$

Then $(K^*)^c \subset K^c$. Recall that $\beta(x) := \beta_1 1_K(x) + \beta_2 1_{K^c}(x) \leq \beta_2$ and let $V \in C^\gamma$ ($\gamma \in (0, 1]$) with

$$\beta_2 1_{(K^*)^c} \leq V \leq \beta. \tag{14}$$

(The existence of a *continuous* function satisfying (14) would of course immediately follow from Uryson's lemma (Lemma 4.4 in [27]). In fact it is easy to see the existence of such functions which are even $C^\infty$ by writing $\beta = \beta_2 - (\beta_2 - \beta_1) 1_K$ and considering the function $V := \beta_2 - (\beta_2 - \beta_1) f$, where $f \geq 1_{K^*}$ and $f$ is a $C^\infty$-function obtained as follows. $f$ is a sum of compactly supported $C^\infty$-functions $f_n$, $n \geq 1$, with disjoint support, where $\text{supp}(f_n)$ is in the $\varepsilon_n$-neighborhood of the $n$th connected component of $1_{K^*}$, with appropriately small $0 < \varepsilon_n$'s.) Consider the operator $L + V$ on $\mathbb{R}^d$ and let $\lambda_c = \lambda_c(\omega)$ denote its generalized principal eigenvalue. Since $V \in C^\gamma$, we are in the setting of Chapter 4 in [32] and Lemma 10. In particular, since $V \leq \beta_2$, one has $\lambda_c \leq \lambda_c(L) + \beta_2$ for every $\omega$ (see Theorem 4.4.4 in [32]).

On the other hand, one gets a lower estimate on $\lambda_c$ as follows. Fix $R > 0$. Since $\beta_2 1_{(K^*)^c} \leq V$, by the homogeneity of the Poisson point process, for almost every environment the set $\{x \in \mathbb{R}^d \mid V(x) = \beta_2\}$ contains a clearing of radius $R$. Hence, by comparison, $\lambda_c \geq \lambda^{(R)}$, where $\lambda^{(R)}$ is the principal Dirichlet eigenvalue of $L + \beta_2$ on a ball of radius $R$. Since $R$ can be chosen arbitrarily large and since $\lim_{R \uparrow \infty} \lambda^{(R)} = \lambda_c(L) + \beta_2$, we conclude that $\lambda_c \geq \lambda_c(L) + \beta_2$ for almost every environment.

---

[4]$B \subset\subset D$ means that $B$ is bounded and its closure $\overline{B}$ satisfies $\overline{B} \subset D$.

From the lower and upper estimates, we obtain that

$$\lambda_c = \lambda_c(L) + \beta_2 \quad \text{for a.e. } \omega. \tag{15}$$

Consider now the branching processes with the same motion component $L$ but with rate $V$, respectively constant rate $\beta_2$. The statements (i) and (ii) of Claim 7 are true for these two processes by (15) and Lemma 10. As far as the original process (with rate $\beta$) is concerned, (i) and (ii) of Claim 7 now follow by $\omega$-wise comparison.

## 6. Proof of Theorem 1

We give an upper and a lower estimate separately.

### 6.1. Upper estimate

Let $\varepsilon > 0$. Using the Markov inequality along with the expectation formula (1), we have that on a set of full **P**-measure:

$$\begin{aligned}
P^\omega[(\log t)^{2/d}(r_t - \beta_2) + c(d,\nu) &> \varepsilon] \\
&= P^\omega\{|Z_t| > \exp[t(\beta_2 - c(d,\nu)(\log t)^{-2/d} + \varepsilon(\log t)^{-2/d})]\} \\
&\leq E^\omega|Z_t| \cdot (\exp[t(\beta_2 - c(d,\nu)(\log t)^{-2/d} + \varepsilon(\log t)^{-2/d})])^{-1} \\
&= \exp[-\varepsilon t(\log t)^{-2/d} + o(t(\log t)^{-2/d})] \to 0 \quad \text{as } t \to \infty.
\end{aligned}$$

### 6.2. Lower estimate

We give a 'bootstrap argument': we first prove a weaker estimate which we will later use to prove the stronger result. For better readability, we broke the relatively long proof into three parts.

#### 6.2.1. Part I: A rough exponential lower estimate
Let $0 < \delta < \beta_2$. Then on a set of full **P**-measure

$$\lim_{t \to \infty} P^\omega(|Z_t| \geq e^{\delta t}) = 1. \tag{16}$$

To prove (16), we invoke the definition of the function $V$ from Section 5.3: $V \in C^\gamma$ ($\gamma \in (0,1]$) with

$$\beta_2 1_{(K^*)^c} \leq V \leq \beta. \tag{17}$$

(Recall that $\beta(x) := \beta_1 1_K(x) + \beta_2 1_{K^c}(x)$.) By comparison, it is enough to prove (16) for the 'smooth version' of the process, where $\beta$ is replaced by $V$. The law of this modified process will be denoted by $P^V$ (and the notation $Z$ is unchanged).

Considering the operator $L + V$ on $\mathbb{R}^d$ we have seen in Section 7 that its generalized principal eigenvalue is $\lambda_c(\Delta; \mathbb{R}^d) + \beta_2 = \beta_2$ for every $\omega$.

Take $R > 0$ large enough so that $\lambda_c = \lambda_c(\frac{1}{2}\Delta + V; B(0,R))$, the principal eigenvalue of $\frac{1}{2}\Delta + V$ on $B(0,R)$ satisfies

$$\lambda_c > \delta.$$

Let $\hat{Z}^R$ be the process obtained from $Z$ by introducing killing at $\partial B(0,R)$ (the corresponding law will be denoted by $P_x^{(R)}$). Then

$$\lim_{t \to \infty} P^V(|Z_t| < e^{\delta t}) \leq \lim_{t \to \infty} P^{(R)}(|\hat{Z}_t^R| < e^{\delta t}). \tag{18}$$

Let $0 \leq \phi = \phi^R$ be the Dirichlet eigenfunction corresponding to $\lambda_c$ on $B(0,R)$, and normalize it by $\sup_{x \in B(0,R)} \phi(x) = 1$. Then we can continue inequality (18) with

$$\leq \lim_{t \to \infty} P^{(R)}(\langle \hat{Z}_t^R, \phi \rangle < e^{\delta t}),$$

where $\langle \hat{Z}_t^R, \phi \rangle := \sum_i \phi(\hat{Z}_t^{R,i})$ and $\{\hat{Z}_t^{R,i}\}$ is the '$i$th particle' in $\hat{Z}_t^R$. Notice that $M_t = M_t^\phi := e^{-\lambda_c t} \langle \hat{Z}_t^R, \phi \rangle$ is a nonnegative martingale (see p. 84 in [15]), and define

$$N := \lim_{t \to \infty} M_t.$$

Since $\lambda_c(B(0,R)) > \delta$, and thus $\lim_{t \to \infty} P^{(R)}(M_t < e^{(\delta - \lambda_c)t} \cap \{N > 0\}) = 0$, the estimate is then continued as

$$= \lim_{t \to \infty} P^{(R)}(M_t < e^{(\delta - \lambda_c)t} | N = 0) P^{(R)}(N = 0) \leq P^{(R)}(N = 0).$$

We have that

$$\lim_{t \to \infty} P^V(|Z_t| < e^{\delta t}) \leq P^{(R)}(N = 0)$$

holds for *all* $R$ large enough. Therefore, in order to prove (16), it is sufficient to show that

$$\lim_{R \to \infty} P^{(R)}(N > 0) = 1. \tag{19}$$

Consider now the elliptic boundary value problem (which of course depends on $K$),

$$\begin{aligned}
\frac{1}{2}\Delta u + V(u - u^2) &= 0 \quad \text{in } B(0,R), \\
\lim_{x \to \partial B(0,R)} u(x) &= 0, \\
u &> 0 \quad \text{in } B(0,R).
\end{aligned} \tag{20}$$

The existence of a solution follows from the fact that $\lambda_c > 0$ by an analytical argument given in [33] pp. 262 and 263. (The proof of the existence of such a $u$ relies on finding so-called lower and upper solutions. The assumption $\lambda_c > 0$ enters the stage when a positive lower solution is constructed.) Uniqueness follows by the semilinear maximum principle (Proposition 7.1 in [16]).

The argument below gives a *probabilistic construction* for the solution. Namely, we show that $w_R(x) := P_x^{(R)}(N > 0)$ solves (20). To see this, let $v = v_R := 1 - w_R$. Let us fix an arbitrary time $t > 0$. Using the Markov and branching properties of $Z$ at time $t$, it is straightforward to show that

$$P^{(R)}(N = 0 | \mathcal{F}_t) = \prod_i P_{\hat{Z}_t^{R,i}}^{(R)}(N = 0).$$

Since the left hand-side of this equation defines a $P^{(R)}$-martingale in $t$, so does the right-hand side. That is

$$\widetilde{M}_t := \prod_i v(\hat{Z}_t^{R,i})$$

defines a martingale. From this, it follows by Theorem 17 of [14] that $v$ solves the equation obtained from the first equation of (20) by switching $u - u^2$ to $u^2 - u$. Consequently, $w_R(x) := P_x^{(R)}(N > 0)$ solves the first equation of (20) itself. That $w_R$ solves the second equation, follows easily from the continuity of Brownian motion. Finally its positivity (the third equation of (20)) follows again from the fact that $\lambda_c > 0$ (see Lemma 6 in [15]).

By the semilinear elliptic maximum principle (Proposition 7.1 in [16]; see also [33]), $w_R(\cdot)$ is monotone increasing in $R$. Using standard arguments, one can show that $0 < w := \lim_{R \to \infty} w_R$ too solves the equation in array (20) (see the proof of Theorem 1 in [33]).

Applying the strong maximum principle to $v := 1 - w$, it follows that $w$ is either one everywhere or less than one everywhere. We now suppose that $0 < w < 1$ and will get a contradiction.

When $d \leq 2$, this is simple. We have

$$\frac{1}{2} \Delta w = V(w^2 - w) \lneqq 0 \quad \text{in } \mathbb{R}^d,$$

($\Delta w$ is nonnegative and not identically zero) and this contradicts the recurrence of the Brownian motion in one and two dimensions, because it is known that if $L$ corresponds to a recurrent diffusion then there is no positive function $w$ satisfying $Lw \lneqq 0$ (see Theorem 4.3.9 in [32]). This contradiction proves that in fact $w = 1$ and consequently it proves (16).

Consider now the problem

$$\frac{1}{2} \Delta u + V(u - u^2) = 0 \quad \text{in } \mathbb{R}^d,$$
$$0 < u < 1.$$
(21)

The nonexistence of solutions in the general $d \geq 1$ case is more subtle than for $d \leq 2$. It follows from the fact that $\beta$ is bounded from below by $\beta_1$ along with Theorem 1.1 and Remark 2.4 in [17] (set $g \equiv \beta_1$ in Remark 2.4 in [17]).

We have now shown (16) and completed the first part of our 'bootstrap' proof.

*6.2.2. Part II: Scaling for the refined lower estimate*

Let us return to the proof of the lower estimate in Theorem 1. Let $\varepsilon > 0$. We have to show that on a set of full **P**-measure,

$$\lim_{t \to \infty} P^\omega[(\log t)^{2/d}(r_t - \beta_2) + c(d, \nu) < -\varepsilon] = 0. \tag{22}$$

To achieve this, we will define a particular function $p_t$ (the definition is given in (31)) satisfying that as $t \to \infty$,

$$p_t = \exp\left[-c(d, \nu) \frac{t}{(\log t)^{2/d}} + \mathrm{o}\left(\frac{t}{(\log t)^{2/d}}\right)\right]. \tag{23}$$

Using this function we are going to show a statement implying (22), namely, that for all $\varepsilon > 0$ there is a set of full **P**-measure, where

$$\lim_{t \to \infty} P^\omega[\log |Z_t| < \beta_2 t + \log p_t - \varepsilon t (\log t)^{-2/d}] = 0. \tag{24}$$

Let us first give an outline of the *strategy of our proof*. A key step will be introducing three different time scales, $\ell(t)$, $m(t)$ and $t$ where $\ell(t) = \mathrm{o}(m(t))$ and $m(t) = \mathrm{o}(t)$ as $t \to \infty$. For the first, shortest time interval, we will use that there are 'many' particles produced and they are not moving 'too far away,' for the second (of length $m(t) - \ell(t)$) we will use that one particle moves into a clearing of a certain size at a certain distance, and in the third one (of length $t - m(t)$) we will use that there is a branching tree emanating from that particle so that a certain proportion of particles of that tree stay in the clearing with probability tending to one.

To carry out this program, first recall the following fact (for example see Lemma 4.5.2 in the proof of Theorem 4.5.1 in [35]): Let

$$R_0 = R_0(d,\nu) := \sqrt{\frac{\lambda_d}{c(d,\nu)}} = \left(\frac{d}{\nu\omega_d}\right)^{1/d}$$

(recall that $\lambda_d$ is the principal Dirichlet eigenvalue of $-\frac{1}{2}\Delta$ on the $d$-dimensional unit ball, and $\omega_d$ is the volume of that ball) and let

$$\rho(l) := R_0 (\log l)^{1/d} - (\log \log l)^2, \quad l > 1. \tag{25}$$

Then,

$$\mathbf{P}(\exists l_0(\omega) > 0 \text{ such that } \forall l > l_0(\omega) \ \exists \text{ clearing } B(x_0, \rho(l)) \text{ with } |x_0| \le l) = 1. \tag{26}$$

Next, let $\ell$ and $m$ be two functions $\mathbb{R}_+ \to \mathbb{R}_+$ satisfying the following:

(i) $\lim_{t\to\infty} \ell(t) = \infty$,
(ii) $\lim_{t\to\infty} \frac{\log t}{\log \ell(t)} = 1$,
(iii) $\ell(t) = o(m(t))$ as $t \to \infty$,
(iv) $m(t) = o(\ell^2(t))$ as $t \to \infty$,
(v) $m(t) = o(t(\log t)^{-2/d})$ as $t \to \infty$.

Note that (i)–(v) are in fact not independent, because (iv) follows from (ii) and (v). For example the following choices of $\ell$ and $m$ satisfy (i)–(v): let $\ell(t)$ and $m(t)$ be arbitrarily defined for $t \in [0, e]$, and

$$\ell(t) := t^{1-1/(\log \log t)}, \qquad m(t) := t^{1-1/(2\log \log t)} \quad \text{for } t \ge t_0 > e.$$

For concreteness, let us fix these choices and note that $\rho$ is monotone increasing for large $t$.

*6.2.3. Part III: Completing the refined lower estimate*
Fix $\delta \in (0, \beta_2)$ and define

$$I(t) := \lfloor \exp(\delta \ell(t)) \rfloor.$$

Let $A_t$ denote the following event:

$$A_t := \{|Z_{\ell(t)}| \ge I(t)\}.$$

By (16) we know that on a set of full $\mathbf{P}$-measure,

$$\lim_{t\to\infty} P^\omega(A_t) = 1. \tag{27}$$

By (27), for $t$ fixed we can work on $A_t \subset \Omega$ and consider $I(t)$ many particles at time $\ell(t)$.

As a next step, we need some control on their spatial position. To achieve this, use Remark 1 to compare BBM's with and without obstacles, and then, use the following result taken from Proposition 2.3 in [12] (the stronger, a.s. result is proved in [29]).

Denote $\mathcal{Z}$ the BBM *without* obstacles (and hence with rate $\beta_2$) starting at the origin with a single particle. Let $R(t) = \bigcup_{s \in [0,t]} \text{supp}(\mathcal{Z}(s))$ denote the range of $\mathcal{Z}$ up to time $t$. Let

$$M(t) = \inf\{r > 0: R(t) \subseteq B(0,r)\} \quad \text{for } d \ge 1, \tag{28}$$

be the radius of the minimal ball containing $R(t)$. Then $M(t)/t$ converges to $\sqrt{2\beta_2}$ in probability as $t \to \infty$.

Going back to the set of $I(t)$ many particles at time $\ell(t)$, (28) yields that even though they are at different locations, still for any $\varepsilon' > 0$, with $P^\omega$-probability tending to one, they are all inside the $(\sqrt{2\beta_2} + \varepsilon')\ell(t)$-ball.

Define

$$\rho(t) = \rho(\ell(t)) := R_0[\log \ell(t)]^{1/d} - [\log \log \ell(t)]^2 \quad \text{for } t \geq t_0 > e^e.$$

Recall (26). With **P**-probability one there is a clearing $B = B(x_0, \rho(t))$ such that $|x_0| \leq \ell(t)$, for all large enough $t > 0$. In the sequel we will assume the 'worst case,' when $|x_0| = \ell(t)$. Indeed, it is easy to see that $|x_0| < \ell(t)$ would help in all the arguments below. (Of course, $x_0$ depends on $t$, but this dependence is suppressed in our notation.) By the previous paragraph, with $P^\omega$-probability tending to one, the distance of $x_0$ from *each* of the $I(t)$ many particles is at most

$$(1 + \sqrt{2\beta_2} + \varepsilon')\ell(t).$$

Now, any such particle moves to $B(x_0, 1)$ in another $m(t) - \ell(t)$ time with probability $q_t$, where (using (iii) and (iv) along with the Gaussian density)

$$q_t = \exp\left(-\frac{[(1+\sqrt{2\beta_2}+\varepsilon')\ell(t)]^2}{2[m(t)-\ell(t)]} + o\left(\frac{[(1+\sqrt{2\beta_2}+\varepsilon')\ell(t)]^2}{2[m(t)-\ell(t)]}\right)\right) \to 0 \quad \text{as } t \to \infty.$$

Let the particle positions at time $\ell(t)$ be $z_1, z_2, \ldots, z_{I(t)}$ and consider the independent system of Brownian particles

$$\{W_{z_i}; \ i = 1, 2, \ldots, I(t)\},$$

where $W_{z_i}(0) = z_i$; $i = 1, 2, \ldots, I(t)$. In other words, $\{W_{z_i}; \ i = 1, 2, \ldots, I(t)\}$ just describes the evolution of the $I(t)$ particles picked at time $\ell(t)$ without respect to their possible further descendants and (using the Markov property) by resetting the clock at time $\ell(t)$.

Let $C_t$ denote the following event:

$$C_t := \{\exists i \in \{1, 2, \ldots, I(t)\}, \ \exists 0 \leq s \leq m(t) - \ell(t) \text{ such that } W_{z_i}(s) \in B(x_0, 1)\}.$$

By the independence of the particles,

$$\limsup_{t\to\infty} P^\omega(C_t^c | A_t) = \limsup_{t\to\infty}(1-q_t)^{I(t)} = \limsup_{t\to\infty}[(1-q_t)^{1/q_t}]^{q_t I(t)}. \tag{29}$$

Since (iii) implies that $\frac{\ell^2(t)}{m(t)} = o(\ell(t))$ as $t \to \infty$ and since (i) is assumed, one has

$$q_t e^{\delta \ell(t)} = \exp\left(-\frac{[\ell(t) + (\sqrt{2\beta_2} + \varepsilon')\ell(t)]^2}{2[m(t)-\ell(t)]} + \delta \ell(t) + o(\ell(t))\right) \to \infty \quad \text{as } t \to \infty.$$

In view of this, (29) implies that $\lim_{t\to\infty} P^\omega(C_t^c | A_t) = 0$. Using this along with (27), it follows that on a set of full **P**-measure,

$$\lim_{t\to\infty} P^\omega(C_t) = 1. \tag{30}$$

Once we know (30), we proceed as follows. Recall that $B = B(x_0, \rho(t))$ and that $\{\mathbb{P}_x; \ x \in \mathbb{R}^d\}$ denote the probabilities corresponding to a single generic Brownian particle $W$ (being different from the $W_{z_i}$ above). Let $\sigma_B^{x_0}$ denote the first exit time from $B$:

$$\sigma_B^{x_0} = \sigma_{B(x_0,\rho(t))}^{x_0} := \inf\{s \geq 0 \mid W_s \notin B\}.$$

Abbreviate $t^* := t - m(t)$ and define

$$p_t := \sup_{x \in B(x_0,1)} \mathbb{P}_x(\sigma_B^{x_0} \geq t^*) = \sup_{x \in B(0,1)} \mathbb{P}_x(\sigma_B \geq t^*), \tag{31}$$

where $\sigma_B := \sigma_B^0$. Recall that the radius of $B$ is

$$\rho(t) = R_0 [\log \ell(t)]^{1/d} - \text{o}([\log \ell(t)]^{1/d}) = \sqrt{\frac{\lambda_d}{c(d,\nu)}} [\log \ell(t)]^{1/d} - \text{o}([\log \ell(t)]^{1/d}),$$

and recall the definition of $\lambda_d$ from Claim 3. Then, as $t \to \infty$,

$$p_t = \exp\left[-\frac{\lambda_d \cdot t^*}{\rho^2(t)} + \text{o}\left(\frac{\lambda_d \cdot t^*}{\rho^2(t)}\right)\right] = \exp\left[-c(d,\nu)\frac{t^*}{[\log \ell(t)]^{2/d}} + \text{o}\left(\frac{t^*}{[\log \ell(t)]^{2/d}}\right)\right]. \tag{32}$$

Using (ii) and (v), it follows that in fact

$$p_t = \exp\left[-c(d,\nu)\frac{t}{(\log t)^{2/d}} + \text{o}\left(\frac{t}{(\log t)^{2/d}}\right)\right]. \tag{33}$$

A little later we will also need the following notation:

$$p_s^t := \sup_{x \in B(x_0,1)} \mathbb{P}_x(\sigma_B^{x_0} \geq s) = \sup_{x \in B(0,1)} \mathbb{P}_x(\sigma_B \geq s), \tag{34}$$

with this notation,

$$p_t = p_{t^*}^t.$$

By slightly changing the notation, let $Z^x$ denote the BBM starting with a single particle at $x \in B$; and let $Z^{x,B}$ denote the BBM starting with a single particle at $x \in B$ and with absorbtion at $\partial B$ (and still branching at the boundary at rate $\beta_2$).

Since branching does not depend on motion, $|Z^{x,B}|$ is a non-spatial Yule's process (and of course it does not depend on $x$) and thus for all $x \in B$,

$$\exists N := \lim_{t \to \infty} e^{-\beta_2 t} |Z_t^{x,B}| > 0 \tag{35}$$

almost surely (see Theorems III.7.1–2 in [1]).

Note that some particles of $Z^x$ may re-enter $B$ after exiting, whereas for $Z^{x,B}$ that may not happen. Thus, by a simple coupling argument, one has that for all $t \geq 0$, the random variable $|Z_t^x(B)|$ is stochastically larger than $|Z_t^{x,B}(B)|$.

Recall that our goal is to show (24), and recall also (31) and (33). In fact, we will prove the following, somewhat stronger version of (24): we will show that if the function $\gamma : [0,\infty) \to [0,\infty)$ satisfies $\lim_{t\to\infty} \gamma_t = 0$, then on a set of full **P**-measure,

$$\lim_{t \to \infty} P^\omega(|Z_t| < \gamma_t \cdot e^{\beta_2 t^*} p_t) = 0. \tag{36}$$

Recalling $t^* = t - m(t)$, and setting

$$\gamma_t := \exp\left(m(t) - \varepsilon \frac{t}{(\log t)^{2/d}}\right) \quad \text{for } t \geq t_0 > e,$$

a simple computation shows that (36) yields (24). Note that this particular $\gamma$ satisfies $\lim_{t\to\infty} \gamma_t = 0$ because of the condition (v) on the function $m$.

By the comparison between $|Z_t^x(B)|$ and $|Z_t^{x,B}(B)|$ (discussed in the paragraph after (35)) along with (30) and the Markov property applied at time $m(t)$, we have that

$$\lim_{t \to \infty} P(|Z_t| < \gamma_t \cdot e^{\beta_2 t^*} p_t) \leq \lim_{t \to \infty} \sup_{x \in B} P(|Z_{t^*}^{x,B}(B)| < \gamma_t \cdot e^{\beta_2 t^*} p_t).$$

Consider now the $J(x,t) := |Z^{x,B}_{t^*}|$ many Brownian paths starting at $x \in B$, which are correlated through common ancestry, and let us denote them by $W_1, \ldots, W_{J(x,t)}$. Let

$$n^x_t := \sum_{i=1}^{J(x,t)} 1_{A_i},$$

where

$$A_i := \{W_i(s) \in B, \ \forall 0 \leq s \leq t\}.$$

Then we have to show that

$$L := \lim_{t \to \infty} \sup_{x \in B} P(n^x_t < \gamma_t \cdot e^{\beta_2 t^*} p_t) = 0. \tag{37}$$

Clearly, for all $x \in B$,

$$L = \lim_{t \to \infty} \sup_{x \in B} P\left(\frac{n^x_t}{Ne^{\beta_2 t^*}} < \frac{\gamma_t p_t}{N}\right)$$

$$\leq \lim_{t \to \infty} \sup_{x \in B} \left[ P\left(\frac{n^x_t}{Ne^{\beta_2 t^*}} < \frac{1}{2} p_t\right) + P(N \leq 2\gamma_t) \right]. \tag{38}$$

Using the fact that $\lim_{t \to \infty} \gamma_t = 0$ and that $N$ is almost surely positive,

$$\lim_{t \to \infty} P(N \leq 2\gamma_t) = 0;$$

hence it is enough to show that

$$\lim_{t \to \infty} \sup_{x \in B} P\left(\frac{n^x_t}{Ne^{\beta_2 t^*}} < \frac{1}{2} p_t\right) = 0. \tag{39}$$

The strategy for the rest of the proof is conditioning on the value of the positive random variable $N$ and then using Chebysev's inequality, for which we will have to carry out some variance calculations. Since the particles are correlated through common ancestry, we will have to handle the distribution of the splitting time of the most common ancestor of two generic particles. Doing so, we will prove a lemma, while some further computations will be deferred to Appendix B.

Let $R$ denote the law of $N$ and define the conditional laws

$$P^y(\cdot) := P(\cdot | N = y), \quad y > 0.$$

Then

$$P\left(\frac{n^x_t}{Ne^{\beta_2 t^*}} < \frac{1}{2} p_t\right) = \int_0^\infty R(dy) \, P^y\left(\frac{n^x_t}{ye^{\beta_2 t^*}} < \frac{1}{2} p_t\right).$$

Define the conditional probabilities

$$\widetilde{P}^y(\cdot) := P^y(\cdot | |Z^B_{t,x}| \geq \mu_t) = P(\cdot | N = y, |Z^B_{t,x}| \geq \mu_t), \quad y > 0,$$

where $\mu_t = \mu_{t,y} := \lfloor \frac{3y}{4} e^{\beta_2 t^*} \rfloor$. Recall that (31) defines $p_t$ by taking supremum over $x$ and that $|Z^{x,B}_t|$ in fact does not depend on $x$. One has

$$P\left(\frac{n^x_t}{Ne^{\beta_2 t^*}} < \frac{1}{2} p_t\right) \leq \int_0^\infty R(dy) \left[ \widetilde{P}^y\left(\frac{n^x_t}{ye^{\beta_2 t^*}} < \frac{1}{2} p_t\right) + P^y\left(e^{-\beta_2 t^*} |Z^{x,B}_t| < \frac{3}{4} y\right) \right]. \tag{40}$$

As far as the second term of the integrand in (40) is concerned, the limit in (35) implies that

$$\lim_{t\to\infty} \int_{\mathbb{R}} R(\mathrm{d}y)\, P^y\!\left(\mathrm{e}^{-\beta_2 t^*}|Z_t^{x,B}| < \frac{3}{4}y\right) = \lim_{t\to\infty} P\!\left(\mathrm{e}^{-\beta_2 t^*}|Z_t^{x,B}| < \frac{3}{4}N\right) = 0.$$

Let us now concentrate on the first term of the integrand in (40). In fact, it is enough to prove that for each fixed $K > 0$,

$$\lim_{t\to\infty} \int_{1/K}^{\infty} R(\mathrm{d}y)\, \widetilde{P}^y\!\left(\frac{n_t^x}{y\mathrm{e}^{\beta_2 t^*}} < \frac{1}{2}p_t\right) = 0. \tag{41}$$

Indeed, once we know (41), we can write

$$\lim_{t\to\infty} \int_0^{\infty} R(\mathrm{d}y)\, \widetilde{P}^y\!\left(\frac{n_t^x}{y\mathrm{e}^{\beta_2 t^*}} < \frac{1}{2}p_t\right)$$
$$\leq \lim_{t\to\infty} \int_{1/K}^{\infty} R(\mathrm{d}y)\, \widetilde{P}^y\!\left(\frac{n_t^x}{y\mathrm{e}^{\beta_2 t^*}} < \frac{1}{2}p_t\right) + R\!\left(\left[0,\frac{1}{K}\right]\right) = R\!\left(\left[0,\frac{1}{K}\right]\right). \tag{42}$$

Since this is true for all $K > 0$, thus letting $K \uparrow \infty$,

$$\lim_{t\to\infty} \int_0^{\infty} R(\mathrm{d}y)\, \widetilde{P}^y\!\left(\frac{n_t^x}{y\mathrm{e}^{\beta_2 t^*}} < \frac{1}{2}p_t\right) = 0.$$

Returning to (41), let us pick randomly $\mu_t$ many points out of the $J(x,t)$ many particles – this is almost surely possible under $\widetilde{P}^y$. (Again, 'randomly' means that the way we pick the particles is independent of their genealogy and their spatial position.) Let us denote the collection of these $\mu_t$ many particles by $M_t$. Define

$$\widehat{n}_t^x := \sum_{i\in M_t} 1_{A_i}.$$

Then, one has

$$\widetilde{P}^y\!\left(\frac{n_t^x}{y\mathrm{e}^{\beta_2 t^*}} < \frac{1}{2}p_t\right) \leq \widetilde{P}^y\!\left(\frac{\widehat{n}_t^x}{y\mathrm{e}^{\beta_2 t^*}} < \frac{1}{2}p_t\right). \tag{43}$$

We are going to use Chebyshev's inequality and therefore we now calculate the variance. Recall that $p_t = \sup_{x\in B(0,1)} \mathbb{P}_x(\sigma_B \geq t^*)$. Using that for $x \in B(0,1)$,

$$\mathbb{P}_x(\sigma_B \geq t) - [\mathbb{P}_x(\sigma_B \geq t)]^2 \leq \mathbb{P}_x(\sigma_B \geq t) \leq \mathbb{P}_x(\sigma_B \geq t^*) \leq p_t,$$

one has

$$\widetilde{\mathrm{Var}}^y(\widehat{n}_t^x) \leq \mu_t p_t + \mu_t(\mu_t - 1)\frac{\sum_{(i,j)\in K(t,x)} \widetilde{\mathrm{cov}}^y(1_{A_i}, 1_{A_j})}{\mu_t(\mu_t - 1)},$$

where $K(t,x) := \{(i,j): i\neq j, 1 \leq i,j \leq \mu_t\}$. Now observe that

$$\frac{\sum_{i,j\in K(t,x)} \widetilde{\mathrm{cov}}^y(1_{A_i}, 1_{A_j})}{\mu_t(\mu_t - 1)} = \mathbf{E}\,\widetilde{\mathrm{cov}}^y(1_{A_i}, 1_{A_j}) = (\mathbf{E}\otimes \widetilde{P}^y)(A_i \cap A_j) - p_t^2,$$

where under $\mathbf{P}$ the pair $(i,j)$ is chosen randomly and uniformly over the $\mu_t(\mu_t - 1)$ many possible pairs.

Let $Q^{t,y}$ and $Q^{(t)}$ denote the distribution of the splitting time of the most recent common ancestor of the $i$th and the $j$th particle under $\widetilde{P}^y$ and under $\widetilde{P}$, respectively. By the Markov property applied at this splitting time $s$, one has

$$(\mathbf{E} \otimes \widetilde{P}^y)(A_i \cap A_j) = p_t \int_{s=0}^{t} \int_B p_{t-s,x}^t \widetilde{p}^{(t)}(0,s,\mathrm{d}x)\, Q^{t,y}(\mathrm{d}s),$$

where

$$\widetilde{p}^{(t)}(0,t,\mathrm{d}x) := \mathbb{P}_0(W_t \in \mathrm{d}x | W_z \in B,\ z \le t).$$

By the Markov property applied at time $s$,

$$p_s^t \int_B p_{t-s,x}^t \widetilde{p}^{(t)}(0,s,\mathrm{d}x) = p_t,$$

and thus

$$(\mathbf{E} \otimes \widetilde{P}^y)(A_i \cap A_j) = p_t \int_{s=0}^{t} \frac{p_t}{p_s^t}\, Q^{t,y}(\mathrm{d}s).$$

Hence

$$\widetilde{\mathrm{Var}}^y(\widehat{n}_t^x) \le \mu_t(p_t - p_t^2) + \mu_t(\mu_t - 1)p_t^2 \cdot (I_t - 1), \tag{44}$$

where

$$I_t := \int_{s=0}^{\infty} [p_s^t]^{-1} Q^{t,y}(\mathrm{d}s).$$

Note that this estimate is uniform in $x$ (see the definition of $p_t$ in (31)). Define also

$$J_t := \int_{s=0}^{\infty} [p_s^t]^{-1} Q^{(t)}(\mathrm{d}s).$$

**Lemma 11.**

$$\lim_{t \to \infty} J_t = 1. \tag{45}$$

The proof of this lemma is deferred to the end of this section.

Once we know (45), we proceed as follows. Using Chebyshev's inequality, one has

$$\widetilde{P}^y\left(\frac{\widehat{n}_t^x}{y e^{\beta_2 t^*}} < \frac{1}{2} p_t\right) \le \widetilde{P}^y\left(|\widehat{n}_t^x - E^y \widehat{n}_t^x| > \frac{1}{4} p_t y e^{\beta_2 t^*}\right) \le 16 \frac{\widetilde{\mathrm{Var}}^y(\widehat{n}_t^x)}{p_t^2 y^2 e^{2\beta_2 t^*}}.$$

By (44), we can continue the estimate by

$$\le 16\left[\frac{\mu_t p_t}{p_t^2 y^2 e^{2\beta_2 t^*}} + \frac{1}{2}\mu_t(\mu_t - 1) \cdot y^{-2} e^{-2\beta_2 t^*} \cdot (I_t - 1)\right].$$

Writing out $\mu_t$, integrating against $R(\mathrm{d}y)$, and using that the lower limit in the integral is $1/K$, one obtains the upper estimate

$$\int_{1/K}^{\infty} R(\mathrm{d}y)\, \widetilde{P}^y\left(\frac{\widehat{n}_t^x}{y e^{\beta_2 t^*}} < \frac{1}{2} p_t\right)$$

$$\le 12 K p_t^{-1} e^{-\beta_2 t^*} + \int_{1/K}^{\infty} R(\mathrm{d}y)\, \frac{1}{2}\mu_t(\mu_t - 1) \cdot y^{-2} e^{-2\beta_2 t^*} \cdot (I_t - 1). \tag{46}$$

(Recall that $I_t$ in fact depends on $y$.) Since $\lim_{t\to\infty} p_t e^{\beta_2 t^*} = \infty$, thus the first term on the right-hand side of (46) tends to zero as $t \to \infty$. Recall now that $\mu_t := \lfloor \frac{3y e^{\beta_2 t^*}}{4} \rfloor$. As far as the second term of (46) is concerned, it is easy to see that it also tends to zero as $t \to \infty$, provided

$$\lim_{t\to\infty} \int_0^\infty R(\mathrm{d}y)(I_t - 1) = 0.$$

But $\int_0^\infty R(\mathrm{d}y)(I_t - 1) = J_t - 1$ and so we are finished by recalling (45). Hence (41) follows. This completes the proof of the lower estimate in Theorem 1.

*6.2.4. Proof of Lemma 11*

Since $J_t \geq 1$, thus it is enough to prove that

$$\limsup_{t\to\infty} J_t \leq 1.$$

For $r > 0$ we denote by $\lambda_r^* := \lambda_c(\frac{1}{2}\Delta, B(0,r))$ the principal eigenvalue of $\frac{1}{2}\Delta$ on $B(0,r)$. Since $\lambda_r^*$ tends to zero as $r \uparrow \infty$ we can pick an $R > 0$ such that $-\lambda_R^* < \beta_2$. Let us fix this $R$ for the rest of the proof.

Let us also fix $t > 0$ for a moment. From the probabilistic representation of the principal eigenvalue (see Chapter 4 in [32]) we conclude the following: for $\hat\varepsilon > 0$ fixed there exists a $T(\hat\varepsilon)$ such that for $s \geq T(\hat\varepsilon)$,

$$\log p_s^t \geq (\lambda_{\rho(t)} - \hat\varepsilon)s.$$

Hence, for $\hat\varepsilon > 0$ small enough ($\hat\varepsilon < -\lambda_R^*$) and for all $t$ satisfying $\lambda_{\rho(t)} \geq \lambda_R^* + \hat\varepsilon$ (recall that $\lim_{t\to\infty} \rho(t) = \infty$) and $s \geq T(\hat\varepsilon, t)$,

$$\log p_s^t \geq \lambda_R^* \cdot s. \tag{47}$$

Note that $T(\hat\varepsilon, t)$ can be chosen uniformly in $t$ because[5] $\lim_{t\to\infty} \rho(t) = \infty$, and so we will simply write $T(\hat\varepsilon)$. Furthermore, clearly, $T(\hat\varepsilon)$ can be chosen in such a way that

$$\lim_{\hat\varepsilon \downarrow 0} T(\hat\varepsilon) = \infty. \tag{48}$$

Depending on $\hat\varepsilon$ let us break the integral into two parts:

$$J_t = \int_{s=0}^{T(\hat\varepsilon)} [p_s^t]^{-1} Q^{(t)}(\mathrm{d}s) + \int_{s=T(\hat\varepsilon)}^{t} [p_s^t]^{-1} Q^{(t)}(\mathrm{d}s) =: J_t^{(1)} + J_t^{(2)}.$$

We are going to control the two terms separately.

*Controlling $J_t^{(1)}$*: We show that

There exists the limit $\lim_{t\to\infty} J_t^{(1)} \leq 1$. \hfill (49)

First, it is easy to check that for all $t > 0$, $Q^{(t)}(\mathrm{d}s)$ is absolutely continuous, i.e. $Q^{(t)}(\mathrm{d}s) = g^{(t)}(s)\,\mathrm{d}s$ with some $g^{(t)} \geq 0$. So

$$J_t^{(1)} = \int_{s=0}^{T(\hat\varepsilon)} [p_s^t]^{-1} Q^{(t)}(\mathrm{d}s) = \int_{s=0}^{T(\hat\varepsilon)} [p_s^t]^{-1} g^{(t)}(s)\,\mathrm{d}s.$$

Evidently, one has $[p_s^t]^{-1} \downarrow$ as $t \to \infty$. Also, since $Q^{(t)}([a,b])$ is monotone non-increasing in $t$ for $0 \leq a \leq b$, therefore $g^{(t)}(\cdot)$ is also monotone non-increasing in $t$. Hence, by monotone convergence,

$$\lim_{t\to\infty} J_t^{(1)} = \int_{s=0}^{T(\hat\varepsilon)} g(s)\,\mathrm{d}s = \lim_{t\to\infty} \int_{s=0}^{T(\hat\varepsilon)} g^{(t)}(s)\,\mathrm{d}s \leq \lim_{t\to\infty} \int_{s=0}^{t} g^{(t)}(s)\,\mathrm{d}s = 1,$$

---

[5]Recall that we picked a version of $\rho(t)$ which is monotone increasing for large $t$'s.

where $g := \lim_{t\to\infty} g^{(t)}$.

*Controlling $J_t^{(2)}$*: Recall that

$$\log p_s^t \geq \lambda_R^* \cdot s, \quad \forall s \geq T(\hat{\varepsilon}). \tag{50}$$

Thus,

$$J_t^{(2)} \leq \int_{T(\hat{\varepsilon})}^t \exp(-\lambda_R^* \cdot s) Q^{(t)}(\mathrm{d}s).$$

We will show that

$$\lim_{\hat{\varepsilon}\downarrow 0} \lim_{t\to\infty} \int_{T(\hat{\varepsilon})}^t \exp(-\lambda_R^* \cdot s) Q^{(t)}(\mathrm{d}s) = \lim_{\hat{\varepsilon}\downarrow 0} \lim_{t\to\infty} \int_{T(\hat{\varepsilon})}^t \mathrm{e}^{-(\lambda_R^* + \beta_2)s} \mathrm{e}^{\beta_2 s} Q^{(t)}(\mathrm{d}s) = 0. \tag{51}$$

Recall that $0 < \beta_2 + \lambda_R^*$. In order to verify (51), we will show that given $t_0 > 0$ there exists some $0 < K = K(t_0)$ with the property that

$$g^{(t)}(s) \leq K s \mathrm{e}^{-\beta_2 s}, \quad \text{for } t > t_0, s \in [t_0, t]. \tag{52}$$

Indeed, it will then follow that

$$\lim_{\hat{\varepsilon}\downarrow 0} \lim_{t\to\infty} \int_{T(\hat{\varepsilon})}^t \exp(-\lambda_R^* \cdot s) Q^{(t)}(\mathrm{d}s) = \lim_{T(\hat{\varepsilon})\to\infty} \lim_{t\to\infty} \int_{T(\hat{\varepsilon})}^t \exp(-\lambda_R^* \cdot s) g^{(t)}(s)(\mathrm{d}s)$$

$$\leq K \lim_{T(\hat{\varepsilon})\to\infty} \int_{T(\hat{\varepsilon})}^\infty s \mathrm{e}^{-(\lambda_R^* + \beta_2)s}(\mathrm{d}s) = 0.$$

Recall that $Q^{(t)}$ corresponds to the conditional law $P(\cdot \mid |Z_{t,x}^B| \geq \mu_t)$. We now claim that we can work with $P(\cdot \mid |Z_{t,x}^B| \geq 2)$ instead of $P(\cdot \mid |Z_{t,x}^B| \geq \mu_t)$. This is because if $Q_0^{(t)}$ corresponds to $P(\cdot \mid |Z_{t,x}^B| \geq 2)$, then an easy computation reveals that for any $\varepsilon > 0$ there exists a $\hat{t}_0 = \hat{t}_0(\varepsilon)$ such that for all $t \geq \hat{t}_0$ and for all $0 \leq a < b$,

$$|Q^{(t)}([a,b]) - Q_0^{(t)}([a,b])| \leq 2(1+\varepsilon) Q_0^{(t)}([a,b]);$$

thus, if

$$Q_0^{(t)}(\mathrm{d}s) \leq L s \mathrm{e}^{-\beta_2 s} \mathrm{d}s \quad \text{on } [t_0, t] \text{ for } t > t_0 \tag{53}$$

holds with some $L > 0$, then also

$$Q^{(t)}(\mathrm{d}s) = g^{(t)}(s)\,\mathrm{d}s \leq K s \mathrm{e}^{-\beta_2 s}\,\mathrm{d}s, \quad t > t_0 \vee \hat{t}_0, s \in [t_0, t] \tag{54}$$

holds with $K := L + 2(1 + \varepsilon)$.

The bound (53) is verified in Appendix B.

It is now easy to finish the proof of (45). To make the dependence on $\hat{\varepsilon}$ clear, let us write $J_t^{(i)} = J_t^{(i)}(\hat{\varepsilon})$, $i = 1, 2$. Then by (49), one has that for *all* $\hat{\varepsilon} > 0$,

$$\limsup_{t\to\infty} J_t \leq 1 + \limsup_{t\to\infty} J_t^{(2)}(\hat{\varepsilon}).$$

Hence, (51) yields

$$\limsup_{t\to\infty} J_t \leq 1 + \lim_{\hat{\varepsilon}\downarrow 0} \limsup_{t\to\infty} J_t^{(2)}(\hat{\varepsilon}) \leq 1,$$

finishing the proof of the lemma. □

## 7. Some additional problems

These problems were deferred to the section after the proofs, because the questions themselves refer to parts of the proofs.

**Problem 12.** The end of the proof for the lower estimate in Theorem 1 is basically a version of the Weak Law of Large Numbers. Using SLLN instead (and making some appropriate changes elsewhere), can one get

$$\liminf_{t \to \infty} (\log t)^{2/d} (r_t - \beta_2) \geq -c(d, \eta) \quad \text{a.s.?}$$

**Problem 13.** The question investigated in this paper was the (local and global) growth rate of the population. The next step can be the following: Once one knows the global population size $|Z_t|$, the model can be rescaled (normalized) by $|Z_t|$, giving a population of fixed weight. In other words, one considers the discrete probability measure valued process

$$\widetilde{Z}_t(\cdot) := \frac{Z_t(\cdot)}{|Z_t|}.$$

Then the question of the *shape* of the population for $Z$ for large times is given by the limiting behavior of the random probability measures $\widetilde{Z}_t$, $t \geq 0$. (Of course, not only the particle mass has to be scaled, but also the spatial scales are interesting – see last paragraph.)

Can one for example locate a *unique dominant branch* for almost every environment, so that the total weight of its complement tends to (as $t \to \infty$) zero?

The motivation for this question comes from our proof of the lower estimate for Theorem 1. It seems conceivable that for large times the 'bulk' of the population will live in a clearing within distance $\ell(t)$ and with radius

$$\rho(t) := R_0 [\log \ell(t)]^{1/d} - [\log \log \ell(t)]^2, \quad t \geq 0,$$

where

$$\lim_{t \to \infty} \ell(t) = \infty \quad \text{and} \quad \lim_{t \to \infty} \frac{\ell(t)}{t} = 0 \quad \text{but} \quad \lim_{t \to \infty} \frac{\log t}{\log \ell(t)} = 1.$$

## 8. Proof of Theorem 2

*8.1. Preparation for the proof*

We first state an abstract lemma as follows.

**Lemma 14.** *Given the probability triple $(\Omega, \mathcal{F}, P)$, let $A_1, A_2, \ldots, A_N \in \mathcal{F}$ be events that are positively correlated in the following sense. If $k \leq N$ and $\{j_1, j_2, \ldots, j_k\} \subseteq \{1, 2, \ldots, N\}$ then $\operatorname{cov}(\mathbf{1}_{A_{j_1} \cap A_{j_2} \cap A_{j_3} \cdots \cap A_{j_{k-1}}}, \mathbf{1}_{A_{j_k}}) \geq 0$. Then*

$$P\left(\bigcap_{i=1}^{N} A_i\right) \geq \prod_{i=1}^{N} P(A_i).$$

**Proof.** We use induction on $N$. Let $N = 2$. Then $P(A_1 \cap A_2) \geq P(A_1) P(A_2)$ is tantamount to $\operatorname{cov}(\mathbf{1}_{A_1}, \mathbf{1}_{A_2}) \geq 0$.

If $N+1$ events are positively correlated in the above sense then any subset of them is positively correlated as well. Given that the statement is true for $N \geq 2$, one has

$$P\left(\bigcap_{i=1}^{N+1} A_i\right) \geq P\left(\bigcap_{i=1}^{N} A_i\right) P(A_{N+1}) \geq P(A_{N+1}) \prod_{i=1}^{N} P(A_i) = \prod_{i=1}^{N+1} P(A_i),$$

and so the statement is true for $N+1$. □

**Corollary 15.** *Consider $Z$, the $(L,\beta)$-branching diffusion where $L$ is a second order differential operator corresponding to the diffusion process $Y$ on $\mathbb{R}^d$ with probabilities $\{\mathbf{Q}_x,\ x \in \mathbb{R}^d\}$, and the branching rate $\beta = \beta(\cdot) \geq 0$ is not identically zero. For $t > 0$ let $N_t$ denote the number of particles at $t$ and let $B_t$ be an open set containing the origin.*

*Let the function $g \colon \mathbb{R}_+ \to \mathbb{N}_+$ be so large that $\lim_{t \to \infty} P(|Z_t| \leq g(t)) = 1$. Then the lower estimate*

$$P_0[\mathrm{supp}(Z_s) \in B_t,\ 0 \leq s \leq t] \geq [\mathbf{Q}_0(Y_s \in B_t,\ 0 \leq s \leq t)]^{g(t)} - \mathrm{o}(1)$$

*holds.*

**Proof.** Let us label the particles in a way that does not depend on their motion. We get $N_t$ (correlated) trajectories of $Y$: $Y^{(i)}$, $1 \leq i \leq N_t$. Denote $A_i := (Y_s^{(i)} \in B_t,\ 0 \leq s \leq t)$. When $N_t < g(t)$, consider some additional (positively correlated) 'imaginary' particles – for example by taking $g(t) - N_t$ extra copies of the first particle. We have

$$P_0[\mathrm{supp}(Z_s) \in B_t,\ 0 \leq s \leq t] = P_0\left(\bigcap_{i=1}^{N_t} A_i\right) \geq P_0\left(\bigcap_{i=1}^{g(t)} A_i \cap \{N_t \leq g(t)\}\right)$$

$$\geq P_0\left(\bigcap_{i=1}^{g(t)} A_i\right) - P_0(N_t > g(t)).$$

It is easy to check that $A_1, A_2, \ldots, A_{g(t)}$ are positively correlated, hence one can continue the lower estimate with

$$\geq \prod_{i=1}^{g(t)} P(A_i) - \mathrm{o}(1) = [\mathbf{Q}_0(Y_s \in B_t,\ 0 \leq s \leq t)]^{g(t)} - \mathrm{o}(1).$$

The proof is finished. □

*8.2. Proof of the theorem*

First, the equation (8) follows from the Taylor expansion $\sqrt{1-x} = 1 - \frac{x}{2} + \mathcal{O}(x^2)$, $x \approx 0$.

For (9), recall that $f(t) := c(d,\nu) \frac{t}{(\log t)^{2/d}}$ and the result on the quenched global population: there are roughly $\exp[\beta_2 t - f(t)]$ particles by time $t$. More precisely, on a set of full $\mathbf{P}$-measure,

$$\lim_{t \to \infty} (\log t)^{2/d} (r_t - \beta_2) = -c(d,\nu) \quad \text{in } P^\omega\text{-probability.} \tag{55}$$

In particular, for all $\varepsilon > 0$, as $t \to \infty$,

$$P^\omega(|Z_t| > \mathrm{e}^{\beta_2 t - f(t) + \varepsilon}) = \mathrm{o}(\log t^{-2/d}). \tag{56}$$

The rest is a straightforward computation. We apply Corollary 15 with

$$g(t) := \lfloor \mathrm{e}^{\beta_2 t - f(t) + \varepsilon} \rfloor.$$

Denote $C_t := \{|Z_t| \leq g(t)\}$. Using that $\lim_{t \to \infty} P^\omega(C_t) = 1$ and that

$$n^2(t) = 2t^2[\beta_2 - c(d,\nu)(\log t)^{-2/d}] \tag{57}$$

along with Corollary 15, it follows that for **P**-almost all $\omega$,

$$P^\omega(A_t) \geq \left(1 - \exp\left[-\frac{n^2(t)}{2t}\right]\right)^{\exp[\beta_2 t - f(t) + \varepsilon]} - \mathrm{o}(1)$$

$$= (1 - \exp[-\beta_2 t + f(t)])^{\exp[\beta_2 t - f(t) + \varepsilon]} - \mathrm{o}(1) \longrightarrow \mathrm{e}^{-\mathrm{e}^\varepsilon}, \quad \text{as } t \to \infty.$$

Consequently, for **P**-almost all $\omega$, $\liminf_{t \to \infty} P^\omega(A_t) > 0$.

## Appendix A. Related population models

In this section we explain how our mathematical setting relates to models in biology.

First, one immediately has the following two biological interpretations in mind:

(i) *Migration with unfertile areas (Population dynamics)*: Population moves in space and reproduces by binary splitting, but randomly located reproduction-suppressing areas modify the growth.
(ii) *Fecundity selection (Genetics)*: Reproduction and mutation takes place. Certain randomly distributed genetic types have low fitness: even though they can be obtained by mutation, they themselves do not reproduce easily, unless mutation transforms them to different genetic types. In genetics this phenomenon is called '*fecundity selection.*' Of course, in this setting 'space' means the space of genetic types rather than physical space.

One question of interest is of course the (local and global) growth rate of the population. Once one knows the global population size, the model can be rescaled (normalized) by the global population size, giving a population of unit mass (somewhat similarly to the fixed size assumption in the Moran model or many other models from theoretical biology) and then the question becomes the *shape* of the population.

In the population dynamics setting this latter question concerns whether or not there is a preferred spatial location for the process to populate. In the genetic setting the question is about the existence of a certain kind of genetic type that is preferred in the long run that lowers the risk of low fecundity caused by mutating into less fit genetics types.

Of course, the *genealogical structure* is a very exciting problem to explore too. For example it seems quite possible that for large times the 'bulk' of the population consists of descendants of a single particle that decided to travel far enough (resp. to mutate many times) in order to be in a less hostile environment (resp. in high fitness genetic type area), where she and her descendants can reproduce freely.

For example, a related phenomenon in marine systems [6] is when hypoxic patches form in estuaries because of stratification of the water. The patches affect different organisms in different ways but are detrimental to some of them. They appear and disappear in an effectively stochastic way. This is an actual system that has some features that correspond to the type of assumptions built into our model.

It appears [22] that a very relevant existing ecological context in which to place our model is the so-called '*source-sink theory.*' The basic idea is that some patches of habitat are good for a species (and growth rate is positive) whereas other patches are poor (and growth rate is smaller, or is zero or negative). Individuals can move between patches randomly or according to more detailed biological rules for behavior.

Another kind of scenario where models such as the proposed one would make sense is in systems that are subject to periodic local disturbances [6]. Those would include forests where trees sometimes fall creating gaps (which have various effects on different species but may harm some) or areas of grass or brush which are subject to occasional fires. Again, the effects may be mixed, but the burned areas can be expected to be less suitable habitats for at least some organisms.

Finally we refer the reader to the excellent monograph [7] for a modern introduction to population models (from the PDE point of view).

## Appendix B. Proof of the bound (53)

In this section we give the proof of the bound (53). In fact we prove a precise formula for the distribution of the splitting time of the most recent common ancestor, which, we believe, is of independent interest. The result and its proof are due to W. Angerer and A. Wakolbinger (personal communication).

For simplicity we set $\beta_2 = 1$; the general case is similar. Let us fix $t > 0$. Then for $0 < u < t$, one has

$$Q_0^{(t)}(s \leq u) = \frac{1 - 2ue^{-u} - e^{-2u} + e^{-t}(2u - 3 + 4e^{-u} - e^{-2u})}{(1 - e^{-t})(1 - e^{-u})^2}; \tag{58}$$

and so the density is

$$f^{(t)}(u) := \frac{dQ_0^t}{dl}(u) = 2\frac{e^{-u}(u - 2 + (u+2)e^{-u}) + e^{-t}(1 - 2ue^{-u} - e^{-2u})}{(1 - e^{-t})(1 - e^{-u})^3},$$

where $l$ denotes Lebesgue measure on $[0, t]$.

**Proof of (58).** Consider the Yule population $Y_t := |Z_{t,x}^B|$ and recall that $Q_0^{(t)}$ corresponds to $P(\cdot | Y_t \geq 2)$. The first observation concerns the Yule *genealogy*. Let us pick a pair of individuals from the Yule population at time $t$, assuming that $Y_t = j$, $j \geq 2$. Denote by $I$ the size of the population *just before* the coalescence time $s$ of the two ancestral lines (where 'before' refers to backward time). That is, let $I := Y_{s+dt}$. Using some formulae from [20], we now show that

$$P(I = i) = \frac{j+1}{j-1} \cdot \frac{2}{(i-1)i} \cdot \frac{i-1}{i+1}. \tag{59}$$

Indeed, the distribution of $I$ equals the conditional distribution of $F$ given $F \leq j$, where $F$ is defined as follows. First, the pure birth process $K = (K_i)$ (with respect to 'Yule-time' $i$) is defined in Section 3.5 of [20], and setting $n = 2$,

$$P(K_{i-1} = 1 | K_i = 2) = \frac{2}{i(i-1)}$$

(see formula (4.10) of the paper). Then the 'hitting time' $F$ is defined (in the same section) by $F := \min\{l : K_l = 2\}$. The formula for the distribution of $F$ (formula (2.3) of the paper) now becomes

$$P(F \leq i) = P(K_i = 2) = \frac{i-1}{i+1}, \quad i \geq 2.$$

Let $i \leq j$. Then

$$\begin{aligned} P(I = i) &= P(F = i | F \leq j) = P(K_{i-1} = 1, K_i = 2 | F \leq j) \\ &= P(K_{i-1} = 1, K_i = 2 | K_j = 2) \\ &= \frac{P(K_{i-1} = 1, K_i = 2)}{P(K_j = 2)}. \end{aligned}$$

From the last three displayed formulae one arrives immediately at (59).

Let us now *embed* the 'Yule time' into *real time*. Since a Yule population stemming from $i$ ancestors has a negative binomial distribution, therefore, using the Markov property at times $u$ and $u + du$, one can decompose

$$P(Y_u = i - 1, Y_{u+du} = i, Y_t = j) = p_1 \cdot p_2 \cdot p_3, \tag{60}$$

where

$$p_1 = e^{-u}(1 - e^{-u})^{i-2},$$

$$p_2 = (i - 1)\,du \quad \text{and}$$

$$p_3 = \binom{j-1}{i-1} e^{-(t-u)i}(1 - e^{-(t-u)})^{j-i}.$$

That is,

$$P(Y_u = i - 1, Y_{u+du} = i, Y_t = j)$$
$$= (i - 1) \binom{j-1}{i-1} e^{-(t-u)i}(1 - e^{-(t-u)})^{j-i} e^{-u}(1 - e^{-u})^{i-2}\,du.$$

Since the pair we have chosen coalesce independently from the rest of the population, the random variables $s$ and $I$ are independent. Using the definition of $s$ first and then the independence remarked in the previous sentence, and finally (59) and (60),

$$P(s \in [u, u + du], Y_{s+dt} = i, Y_t = j) = P(I = i, Y_u = i - 1, Y_{u+du} = i, Y_t = j)$$
$$= P(I = i)P(Y_u = i - 1, Y_{u+du} = i, Y_t = j)$$
$$= \binom{j-2}{i-2}\frac{2(j+1)}{i(i+1)} e^{-(t-u)i}(1 - e^{-(t-u)})^{j-i} e^{-u}(1 - e^{-u})^{i-2}\,du,$$

for $0 < u < t$.

Now, summing from $j = i$ to $\infty$, and from $i = 2$ to $\infty$, and then dividing by $P(Y_t \geq 2) = 1 - e^{-t}$, one obtains (after doing some algebra) that for $0 < u < t$,

$$Q_0^{(t)}(s \in (u, u + du)) = \sum_{i=2}^{\infty} e^{-u} \frac{2(2e^{-(t-u)} + i - 1)(1 - e^{-u})^{i-2}}{(1 - e^{-t})i(i+1)}\,du$$
$$= 2 \cdot \frac{e^{-u}(u - 2 + (u + 2)e^{-u}) + e^{-t}(1 - 2ue^{-u} - e^{-2u})}{(1 - e^{-t})(1 - e^{-u})^3}\,du.$$

Equivalently, in integrated form, one has (58). □

## Acknowledgement

I owe thanks to S. C. Harris, A.-S. Sznitman and A. Wakolbinger for helpful discussions, and to two referees for their close reading of the manuscript. As far as mathematical biology is concerned, I am grateful to C. Cosner and B. Fagan for bringing some interesting population models to my attention.

## References

[1] K. B. Athreya and P. E. Ney. *Branching Processes*. Springer, New York, 1972. (Reprinted by Dover, 2004.) MR0373040
[2] S. Albeverio and L. V. Bogachev. Branching random walk in a catalytic medium. I. Basic equations. *Positivity* **4** (2000) 41–100. MR1740207
[3] M. van den Berg, E. Bolthausen and F. den Hollander. Brownian survival among Poissonian traps with random shapes at critical intensity. *Probab. Theory Related Fields* **132** (2005) 163–202. MR2199290
[4] M. D. Bramson. Minimal displacement of branching random walk. *Z. Wahrsch. Verw. Gebiete* **45** (1978) 89–108. MR0510529
[5] M. D. Bramson. Maximal displacement of branching Brownian motion. *Comm. Pure Appl. Math.* **31** (1978) 531–581. MR0494541
[6] C. Cosner. Personal communication.

[7] R. S. Cantrell and C. Cosner. *Spatial Ecology via Reaction-Diffusion Equations*. Wiley, Chichester, 2003. MR2191264
[8] D. Dawson and K. Fleischmann. Catalytic and mutually catalytic super-Brownian motions. In *Proceedings of the Ascona '99 Seminar on Stochastic Analysis, Random Fields and Applications*. R. C. Dalang, M. Mozzi and F. Russo (Eds) 89–110. Birkhäuser, Boston, 2002. MR1958811
[9] M. Donsker and S. R. S. Varadhan. Asymptotics for the Wiener sausage. *Comm. Pure Appl. Math.* **28** (1975) 525–565. MR0397901
[10] E. B. Dynkin. Branching particle systems and superprocesses. *Ann. Probab.* **19** (1991) 1157–1194. MR1112411
[11] J. Engländer. On the volume of the supercritical super-Brownian sausage conditioned on survival. *Stochastic Process. Appl.* **88** (2000) 225–243. MR1767846
[12] J. Engländer and F. den Hollander. Survival asymptotics for branching Brownian motion in a Poissonian trap field. *Markov Process. Related Fields* **9** (2003) 363–389. MR2028219
[13] J. Engländer and K. Fleischmann. Extinction properties of super-Brownian motions with additional spatially dependent mass production. *Stochastic Process. Appl.* **88** (2000) 37–58. MR1761691
[14] J. Engländer and A. E. Kyprianou. Markov branching diffusions: martingales, Girsanov-type theorems and applications to the long term behaviour. Technical report. 1206, Department of Mathematics, Utrecht University, 2001, 39 pages. Available at http://www.math.uu.nl/publications.
[15] J. Engländer and A. E. Kyprianou. Local extinction versus local exponential growth for spatial branching processes. *Ann. Probab.* **32** (2004) 78–99. MR2040776
[16] J. Engländer and R. Pinsky. On the construction and support properties of measure-valued diffusions on $D \subset R^d$ with spatially dependent branching. *Ann. Probab.* **27** (1999) 684–730. MR1698955
[17] J. Engländer and P. L. Simon. Nonexistence of solutions to KPP-type equations of dimension greater than or equal to one. *Electron. J. Differential Equations* (2006) 6 pp. (electronic). MR2198922
[18] J. Engländer and D. Turaev. A scaling limit theorem for a class of superdiffusions. *Ann. Probab.* **30** (2002) 683–722. MR1905855
[19] J. Engländer and A. Winter. Law of large numbers for a class of superdiffusions. *Ann. Inst. H. Poincaré Probab. Statist.* **42** (2006) 171–185. MR2199796
[20] A. Etheridge, P. Pfaffelhuber and A. Wakolbinger. An approximate sampling formula under genetic hitchhiking. *Ann. Appl. Probab.* **16** (2006) 685–729. MR2244430
[21] S. N. Evans and D. Steinsaltz. Damage segregation at fissioning may increase growth rates: a superprocess model. Preprint.
[22] B. Fagan. Personal communication.
[23] K. Fleischmann, C. Mueller and P. Vogt. On the large scale behavior of super-Brownian motion in three dimensions with a single point source. Preprint. Avaible at ArXiv:math.PR/0607670.
[24] M. Freidlin. *Functional Integration and Partial Differential Equations*. Princeton University Press, 1985. MR0833742
[25] J. W. Harris, S. C. Harris and A. E. Kyprianou. Further probabilistic analysis of the Fisher–Kolmogorov–Petrovskii–Piscounov equation: one sided travelling-waves. *Ann. Inst. H. Poincaré Probab. Statist.* **42** (2006) 125–145. MR2196975
[26] H. Kesten and V. Sidoravicius. Branching random walk with catalysts. *Electron. J. Probab.* **8** (2003) (electronic). MR1961167
[27] J. L. Kelley. General topology. *Graduate Texts in Mathematics*. Springer, New York–Berlin, 1975. Reprint of the 1955 edition [Van Nostrand, Toronto, Ont.]. MR0370454
[28] A. Klenke. A review on spatial catalytic branching. *Stochastic Models (Ottawa, ON, 1998)* 245–263. Amer. Math. Soc., Providence, RI, 2000. MR1765014
[29] A. E. Kyprianou. Asymptotic radial speed of the support of supercritical branching and super-Brownian motion in $R^d$. *Markov Process. Related Fields* **11** (2005) 145–156. MR2099406
[30] F. Merkl and M. V. Wüthrich. Infinite volume asymptotics of the ground state energy in a scaled Poissonian potential. *Ann. Inst. H. Poincaré Probab. Statist.* **38** (2002) 253–284. MR1899454
[31] T. Y. Lee and F. Torcaso. Wave propagation in a lattice KPP equation in random media. *Ann. Probab.* **26** (1998) 1179–1197. MR1634418
[32] R. G. Pinsky. *Positive Harmonic Functions and Diffusion*. Cambridge Univ. Press, 1995. MR1326606
[33] R. G. Pinsky. Transience, recurrence and local extinction properties of the support for supercritical finite measure-valued diffusions. *Ann. Probab.* **24** (1996) 237–267. MR1387634
[34] S. Sethuraman. Conditional survival distributions of Brownian trajectories in a one dimensional Poissonian environment. *Stochastic Process. Appl.* **103** (2003) 169–209. MR1950763
[35] A. Sznitman. *Brownian Motion, Obstacles and Random Media*. Springer, Berlin, 1998. MR1717054
[36] J. Xin. Front propagation in heterogeneous media. *SIAM Rev.* **42** (2000) 161–230 (electronic). MR1778352